\documentclass[11pt,reqno]{amsart}
\usepackage[utf8]{inputenc}
\usepackage[T1]{fontenc}
\usepackage[english]{babel}
\usepackage[usenames,dvipsnames]{xcolor}

\usepackage{hyperref}
\definecolor{darkred}{rgb}{0.7,0,0}
\definecolor{darkblue}{rgb}{0,0,0.7}
\definecolor{lightgray}{rgb}{0.96,0.96,1}
\definecolor{commentcolor}{rgb}{1,0,1}

\hypersetup{colorlinks=true,urlcolor=darkred,linkcolor=darkred,citecolor=darkblue}
\usepackage{courier}
\usepackage{listings}
\usepackage[justification=centering]{caption}

\usepackage{dsfont}
\usepackage{amsfonts}
\usepackage{amssymb}
\usepackage{amsmath}
\usepackage{amsthm}
\usepackage[left=3cm, right=3cm, top=2.1cm,bottom=2.1cm]{geometry}
\usepackage{mathtools}
\usepackage{etoolbox}

\patchcmd{\section}{\scshape}{\bfseries}{}{}
\makeatletter
\renewcommand{\@secnumfont}{\bfseries}
\makeatother

\usepackage{cleveref}
\crefname{thm}{Theorem}{Theorems}
\crefname{prop}{Proposition}{Propositions}
\crefname{thmintro}{Theorem}{Theorems}
\crefname{propintro}{Proposition}{Propositions}
\crefname{lemma}{Lemma}{Lemmas}
\crefname{rem}{Remark}{Remarks}
\crefname{cor}{Corollary}{Corollaries}
\crefname{defi}{Definition}{Definitions}
\crefname{ex}{Example}{Examples}
\crefname{section}{Section}{Sections}

\newtheoremstyle{standard}{2ex}{2ex}{\itshape}{}{\bfseries}{.}{.5em}{}
\theoremstyle{standard}
\newtheorem{lemma}{Lemma}[section]
\newtheorem{prop}[lemma]{Proposition}
\newtheorem{cor}[lemma]{Corollary} 
\newtheorem{thm}[lemma]{Theorem}
\newtheorem{thmintro}{Theorem}[section]

\newtheoremstyle{definition}{2ex}{2ex}{}{}{\bfseries}{.}{.5em}{}    
\theoremstyle{definition}
\newtheorem{defi}[lemma]{Definition}
\newtheorem{rem}[lemma]{Remark} 
\newtheorem{ex}[lemma]{Example}

\newtheoremstyle{intermediate}{2ex}{2ex}{}{}{\itshape}{.}{.5em}{}    
\theoremstyle{intermediate}

\DeclareMathOperator{\End}{End}
\DeclareMathOperator{\lcm}{lcm}
\DeclareMathOperator{\red}{{\to}}
\DeclareMathOperator{\ord}{ord}

\newcommand{\IN}{\mathds{N}}
\newcommand{\IZ}{\mathds{Z}}
\newcommand{\IF}{\mathds{F}}
\newcommand{\T}{\mathcal{T}}
\newcommand{\IP}{\mathds{P}}

\title{Equational proofs of Jacobson's Theorem}
\author{Martin Brandenburg}
\thanks{\emph{E-mail address:} \ttfamily{brandenburg@uni-muenster.de}}
\date{\today}

\begin{document}

% Abstract

\begin{abstract}
\noindent A classical theorem by Jacobson says that a ring in which every element $x$ satisfies the equation $x^n=x$ for some $n>1$ is commutative. According to Birkhoff's Completeness Theorem, if $n$ is fixed, there must be an equational proof of this theorem. But equational proofs have only appeared for some values of $n$ so far. This paper is about finding such a proof in general. We are able to make a reduction to the case that $n$ is a prime power $p^k$ and the ring has characteristic $p$. We then prove the special cases $k=1$ and $k=2$. The general case is reduced to a series of constructive Wedderburn Theorems, which we can prove in many special cases. Several examples of equational proofs are discussed in detail.
\end{abstract}
 
\maketitle

% Introduction

\section{Introduction}

Jacobson has shown in \cite{Ja45} that every ring in which every element satisfies an equation of the form
\[x^n = x\]
is commutative, where $n > 1$ may depend on $x$. This theorem has been strengthened by various authors. For example, Herstein \cite{He53} showed that it is enough to require that $x^n - x$ is central. Recently, Anderson and Danchev \cite{AD20} proved the commutativity of rings satisfying the equation $x^n = x^m$, where $n,m$ are of opposite parity.

In the present paper, we are interested in the special case of rings satisfying $x^n = x$ for a fixed number $n > 1$, which we call \emph{$n$-rings}. The case $n=2$ corresponds to the well-studied boolean rings, for which there is a simple equational commutativity proof (see \cref{2case}). More generally, by Birkhoff's Completeness Theorem \cite{Bi35} (see also \cite{Ta79} for a modern account), there must be an equational proof of the commutativity of $n$-rings, which thus exists of simple symbolic manipulations and substitutions that derive $xy=yx$ from the ring axioms and $x^n=x$.

However, Jacobson's proof is not equational, and indeed quite sophisticated. More simple proofs of Jacobson's Theorem have been published by Forsythe and McCoy \cite{FM46}, Herstein \cite{He54,He61}, Rogers \cite{Ro71} and Dolan \cite{Do76}, but arguably the most elementary and slick proofs have been given by Wamsley \cite{Wa71} and Nagahara and Tominaga \cite{NT74} (their proofs are the same). Still, these proofs are not equational.
  
Equational proofs have only been found for some values of $n$ so far. Moreover, they can be quite challenging to find. Buckley and MacHale \cite{BM13} present the classical case $n=2$, several elementary proofs in the case $n=3$, and also some useful lemmas and generalizations. Wavrik \cite{Wa99} discusses the cases $n=3,4$, several similar commutativity problems, and also presents a general algorithm which is supposed to derive $xy=xy$ from $x^n=x$ for general $n$, but the author does not verify that it always works. Zhang \cite{Zh90} presents an algorithm for even $n$, which, in numerical experiments, often reduces the problem to the cases $n=2$ and $n=4$. (The present paper shows that this approach is too optimistic, though.)

Morita's impressive work \cite{Mo78} gives equational proofs for even numbers $\leq 50$ (except for $16,22,32,46$, which are much harder) and odd numbers $\leq 25$. Some larger numbers are covered as well. MacHale \cite{Ma86} covers the case that $n$ is the sum of two distinct proper powers of $2$, even under the weaker assumption that $x^n-x$ is central. Buckley and MacHale have continued this study in \cite{BM12}. Forsythe and McCoy \cite{FM46} give an equational proof that $p$-rings of characteristic $p$ are commutative, where $p$ is a prime number. The general case, however, seems to be an open problem, and no systematic approach has been taken so far.
  
If we follow Birkhoff's proof in our case, we get the following evidence for the existence of an equational proof: Consider the free ring in two variables $\IZ\langle X,Y\rangle$, consisting of non-commuting polynomials in $X$ and $Y$. Then the quotient ring $\IZ\langle X,Y\rangle / \langle f^n - f : f \in \IZ \langle X,Y \rangle\rangle$ is an $n$-ring. By Jacobson's Theorem, this ring is commutative. In particular, the image of the commutator $XY-YX$ vanishes. Hence, there must be an equation
\begin{equation}\label{birkhoffliko}
XY - YX = \sum_{i=1}^{s} g_i \cdot (f_i^n -f_i) \cdot h_i
\end{equation}
for $f_i,g_i,h_i \in \IZ\langle X,Y \rangle$ for $i=1,\dotsc,s$. Once it has been found, this non-commutative polynomial equation can be verified from the ring axioms alone. Evaluating this equation in any $n$-ring yields an equational proof that it is commutative, since the right hand side vanishes.

However, that approach does not tell us at all how $f_i,g_i,h_i$ look like explicitly, and the existing proofs of Jacobson's Theorem do not \emph{construct} $f_i,g_i,h_i$. In fact, the proofs use non-constructive methods such as the axiom of choice (with the exception of \cite{Wa71,NT74,Do76}) and the law of the excluded middle (LEM).

The LEM states that for every statement $P$, either $P$ or $\neg P$ holds. Equivalently, $\neg \neg P$ implies $P$. For example, notice that \cite{Wa71} merely shows $\neg \neg P$ when $P$ is the statement that $n$-rings are commutative. This also means that Jacobson's Theorem has not been proven so far in intuitionistic logic, which is a type of constructive logic that does not include LEM \cite{Da01}.

For an equational proof, only constructive methods are allowed, and we need to \emph{find} an equation like \eqref{birkhoffliko}. On the other hand, such a single equation is quite complicated in practice, and an equational proof should better be divided into multiple steps that are easier to motivate and to digest. For example, the equation
\begin{align*}
XY - YX  & = \bigl((X+Y)^2 - (X+Y)\bigr) - \bigl(X^2-X\bigr) - \bigl(Y^2-Y\bigr) \\
& \quad + \bigl((YX)^2 - (YX)\bigr) - \bigl((-YX)^2 - (-YX)\bigr),
\end{align*}
in $\IZ\langle X,Y \rangle$ proves that every $2$-ring is commutative. But it is better to present the proof in a more structured way (as in \cref{2case}). The case of $3$-rings illustrates this even better (see \cref{3caselong}).

Our goal is to develop a method that works for every value of $n$. We also prove several reduction results, which generalize observations from \cite{Mo78,Ma86,Zh90}. They are motivated by the well-known classification of commutative $n$-rings, which we recall in \cref{sec:struct}.
 
Our first result is a reduction to the case of prime characteristic. For reasons to be explained in \cref{charmeaning}, we write $p=0$ instead, which also includes the zero ring.

\begin{thmintro}[\cref{proofdecomp}]
Let $n > 1$ and assume that for every prime number $p$ with $p-1 \mid n-1$ there is an equational proof that every $n$-ring with $p=0$ is commutative. Then there is an equational proof that every $n$-ring is commutative.
\end{thmintro}

The reader might object that this theorem is true since Jacobson's Theorem in conjunction with Birkhoff's Theorem shows that there is an equational proof anyway. But we have already observed that this approach does not \emph{construct} a proof. The actual meaning of our theorem is that its proof is a method that, given equational proofs with the additional assumption $p=0$, outputs an equational proof of the general case. The goal of the theorem is not to verify that a certain statement is true in the classical sense. Its goal is to construct equational proofs. And of course, our theorem does not assume Jacobson's Theorem to be true.
  
The next result concerns a useful reduction to prime powers.
 
\begin{thmintro}[\cref{reduc}]
Let $n>1$ and let $p$ be a prime number with $p-1 \mid n-1$. Let $k$ be the $\lcm$ of all $d \geq 1$ such that $p^{d}-1 \mid n-1$. Then there is an equational proof that every $n$-ring with $p=0$ is a $p^k$-ring.
\end{thmintro}

Hence, in order to give an equational proof that $n$-rings are commutative, one may restrict to $p^k$-rings with $p=0$. Again, the actual content of the theorem is the method that results from its proof, which we then apply in several examples. For example, every $9$-ring with $5=0$ is a $5$-ring, and every $10$-ring is a $4$-ring. It is not easy to prove this directly without any guidance, but our proof is the outcome of a general method, so that no individual ideas are necessary to find it.

Morita \cite{Mo78} proved several instances of \cref{reduc}: Direct and clever calculations show that $n$-rings coincide with $2$-rings for $n = 6,12,14,18,20,24,26,30,38,42,44,48$, that $n$-rings coincide with $4$-rings for $n=10,28,34,40$, that $n$-rings coincide $8$-rings for $n=36,50$, and that $46$-rings coincide with $16$-rings. (Some of these calculations are omitted, though.) In contrast, our approach does not depend on a specific value of $n$. Incidentally, the classification of the numbers $n$ such that $n$-rings with $2=0$ coincide with $2$-rings can be found in \cite{HLY94}.
 
The next result concerns the special case $k=1$.
 
\begin{thmintro}[\cref{pcase}]
Let $p$ be a prime number. If $A$ is $p$-ring with $p=0$, then every element of $A$ can be explicitly written as a $\IZ$-linear combination of idempotent elements. In particular, $A$ is commutative by an equational proof.
\end{thmintro}

The last conclusion comes from the basic fact that in a reduced ring, idempotent elements are central. We also provide a version of that theorem that works over any finite field.

Sometimes the case $k=1$ is enough: Let us call a natural number $n > 1$ \emph{simple} if every prime power $q$ with $q-1 \mid n-1$ is actually a prime number. There are many simple numbers, the first ones are $2, 3, 5, 6, 11, 12, 14, 18, 20$, and for these, Jacobson's Theorem is easy to prove (\cref{jacsimple}). We can turn it into an equational proof as well:
 
\begin{thmintro}[\cref{speceq}] 
For every simple number $n>1$ there is an equational proof that every $n$-ring is commutative.
\end{thmintro}

Even though the recipe is always the same, these equational proofs are not ``uniform'' in $n$, though. They have to be written down for every single value of $n$, and for large values of $n$ this can be quite cumbersome. We do not know if there is a uniform proof.

For non-simple numbers, our method is often able to reduce the problem to smaller numbers. For example, the commutativity of $2023$-rings can be reduced to the commutativity of $4$-rings (see \cref{2023case}).

We also give a proof of the case $k=2$, inspired by Morita's proofs for $n=9,25$.

\begin{thmintro}[\cref{p2case}]
Let $p$ be a prime number. There is an equational proof that every $p^2$-ring with $p=0$ is commutative.
\end{thmintro}

In conjunction with the previous reduction results, this covers many more examples. For example, we find an equational proof that $73$-rings are commutative (\cref{73case}).

As for the general case, we take inspiration from the elementary proof of Jacobson's Theorem in \cite{Wa71,NT74}. It reduces the problem to finite rings, where Wedderburn's Theorem comes into play, which states that finite division rings are commutative. We are able to make the reduction constructive and even equational. We just have to replace Wedderburn's Theorem by the following constructive variant: If $f \in \IF_p[T]$, then $W_{p,k,f}$ states that in a $p^k$-ring with $p=0$, we have
\[ba=f(a)b \implies ba=ab.\]
In classical logic, Wedderburn's Theorem implies $W_{p,k,f}$ (\cref{Wc}). We then prove:

\begin{thmintro}[\cref{generalcase-eq,weddermonomial}]
Let $p$ be a prime number and $k \geq 1$. Assume that for all $1 \leq m < k$ there is an equational proof of the constructive Wedderburn Theorem $\smash{W_{p,k,T^{p^m}}}$. Then there is an equational proof that every $p^k$-ring with $p=0$ is commutative. The same conclusion holds when there are equational proofs of all $W_{p,k,f}$ where $f$ is of degree $<k$.
\end{thmintro}

Thus, the problem reduces to finding proofs of certain constructive Wedderburn Theorems. We do not succeed in general, but the following special case already covers many cases.

\begin{thmintro}[\cref{gcdmain}]
Let $p$ be a prime number and $k \geq 1$. If $\gcd(k,p^k-1)=1$, then there is an equational proof that every $p^k$-ring with $p=0$ is commutative.
\end{thmintro}
  
Table \ref{resulttable} summarizes our results for ${2 \leq n \leq 101}$. Writing $\red m_1,m_2,\dotsc$ means that commutativity can be reduced to the classes of $m_1$-rings, $m_2$-rings, etc.\ (by \cref{proofdecomp,reduc,pcase}). Writing $=\!m$ means that the classes of $n$-rings and $m$-rings coincide. An empty cell means that the equational proof is not complete yet. We plan to fill these gaps in future work.

\begin{center}
\begin{table}
\begin{tabular}[c]{ll|ll|ll|ll|ll}
$n$ & proof & $n$ & proof & $n$ & proof & $n$ & proof & $n$ & proof \\ \hline

$2$ & \eqref{speceq} &
$22$ & $\red 64$ &
$42$ & $=2$ &
$62$ & $=2$ &
$82$ & $=4$ \\

$3$ & \eqref{speceq} &
$23$ & \eqref{speceq} &
$43$ & $\red 64$ &
$63$ & $\red 32$ &
$83$ & \eqref{speceq} \\

$4$ & \eqref{p2case} &
$24$ & $=2$ &
$44$ & $=2$ &
$\mathbf{64}$ & &
$84$ & $=2$ \\

$5$ & \eqref{speceq} &
$25$ & \eqref{p2case} &
$45$ & \eqref{speceq} &
$65$ & $\red 9$ &
$85$ & $\red 64$ \\

$6$ & $=2$ &
$26$ & $=2$ &
$46$ & $=16$ &
$66$ & $=2$ &
$86$ & $=2$ \\

$7$ & $\red 4$ &
$27$ & \eqref{gcdmain} &
$47$ & \eqref{speceq} &
$67$ & $\red 4$ &
$87$ & $=3$ \\

$8$ & \eqref{8case} &
$28$ & $=4$ &
$48$ & $=2$ &
$68$ & $=2$ &
$88$ & $=4$ \\

$9$ &\eqref{p2case} &
$29$ & $\red 8$ &
$49$ & \eqref{p2case} &
$69$ & \eqref{speceq} &
$89$ & $\red 9$ \\

$10$ & $=4$ &
$30$ & $=2$ &
$50$ & $=8$ &
$70$ & $=4$ &
$90$ & $=2$ \\

$11$ & \eqref{speceq} &
$31$ & $\red 16$ &
$51$ & \eqref{speceq} &
$71$ & $\red 8$ &
$91$ & $\red 16$ \\

$12$ & $=2$ &
$32$ & \eqref{gcdmain} &
$52$ & $=4$ &
$72$ & $=2$ &
$92$ & $=8$ \\

$13$ & $\red 4$ &
$33$ & $\red 9$ &
$53$ & $\red 27$ &
$73$ & $\red 4,9,25$ &
$93$ & \eqref{speceq} \\

$14$ & $=2$ &
$34$ & $=4$ &
$54$ & $=2$ &
$74$ & $=2$ &
$94$ & $\red 1024$ \\

$15$ & $\red 8$ &
$35$ & $=3$ &
$55$ & $\red 4$ &
$75$ & \eqref{speceq} &
$95$ & $=3$ \\

$16$ & \eqref{gcdmain} &
$36$ & $=8$ &
$56$ & $=2$ &
$76$ & $=16$ &
$96$ & $=2$ \\

$17$ & $\red 9$ &
$37$ & $\red 4$ &
$57$ & $\red 8,9$ &
$77$ & \eqref{speceq} &
$97$ & $\red 4,9,25,49$ \\

$18$ & $=2$ &
$38$ & $=2$ &
$58$ & $=4$ &
$78$ & $=8$ &
$98$ & $=2$ \\

$19$ & $\red 4$ &
$39$ & $=3$ &
$59$ & \eqref{speceq} &
$79$ & $\red 4,27$ &
$99$ & $\red 8$ \\

$20$ & $=2$ &
$40$ & $=4$ &
$60$ & $=2$ &
$80$ & $=2$ &
$100$ & $=4$ \\

$21$ & \eqref{speceq} &
$41$ & $\red 9$ &
$61$ & $\red 16$ &
$\mathbf{81}$ & &
$101$ & \eqref{speceq}

\end{tabular}
\caption{Equational proofs for $2 \leq n \leq 101$ and their reductions}\label{resulttable}
\end{table}
\end{center}
\vspace{-5mm}

% Table of contents

\tableofcontents

% Preliminaries

\section{Preliminaries}
 
Following Jacobson \cite[Section 2.17]{Ja12}, rings are assumed to be unital, and we explicitly speak of rings without a required unit when necessary, which we call \emph{rngs}. See also \cite{Po14}. Of course, rings are not assumed to be commutative unless otherwise stated. The set of prime numbers is denoted by $\IP$.

We need to briefly formalize the notion of an equational proof in order to clarify what is allowed in our proofs and what not. For more details, we refer to \cite{Bi35} and \cite{Ta79}.

\begin{defi} \label{eqlogic}
If $\T$ is an algebraic theory and $\sigma,\tau$ are terms in $\T$, by an \emph{equational proof of $\sigma = \tau$} we mean a deduction $\T \vdash \sigma = \tau$ based on the following rules:
\begin{enumerate}
\item\label{rule1} If $\sigma = \tau$ is an axiom of $\T$, then $\T \vdash \sigma = \tau$.
\item\label{rule2} We have $\T \vdash \sigma = \sigma$.
\item\label{rule3} If $\T \vdash \sigma = \tau$, then $\T \vdash \tau = \sigma$.
\item\label{rule4} If $\T \vdash \sigma = \tau$ and $\T \vdash \tau = \rho$, then $\T \vdash \sigma = \rho$.
\item\label{rule5} If $\rho$ is a term of arity $n$ in $\T$ and $\sigma_1,\dotsc,\sigma_n$, $\tau_1,\dotsc,\tau_n$ are terms in $\T$ with $\T \vdash \sigma_i = \tau_i$ for all $i=1,\dotsc,n$, then $\T \vdash \rho(\sigma_1,\dotsc,\sigma_n) = \rho(\tau_1,\dotsc,\tau_n)$.
\item\label{rule6} If $\sigma,\tau$ are terms of arity $n$ in $\T$ with $\T \vdash \sigma = \tau$, then $\T \vdash \sigma(\rho_1,\dotsc,\rho_n) = \tau(\rho_1,\dotsc,\rho_n)$ for all terms $\rho_1,\dotsc,\rho_n$ in $\T$.
\end{enumerate}
Rule \eqref{rule1} is self-explanatory. Rules \eqref{rule2},\eqref{rule3},\eqref{rule4} state that equational equality is an equivalence relation. Rule \eqref{rule5} means that we may replace equationally equal terms in expressions. Rule \eqref{rule6} allows substitutions of variables with more complex expressions. In practice, we will just write $\sigma = \tau$ instead of $\T \vdash \sigma = \tau$. Notice that there are no quantors and no negation.
\end{defi}

\begin{ex}
There is an equational proof that every group $G$ with $g^2=1$ for all $g \in G$ is commutative. More precisely, in the theory of groups with the additional axiom $g^2=1$, we derive $gh=hg$ as follows: We have $(gh)^2=1$ by rule \eqref{rule6}, and $1=1 \cdot 1 = g^2 \cdot h^2$ by rule \eqref{rule5}. (Of course, we have also applied the other, more trivial rules.) Hence, $(gh)^2=g^2 h^2$, which means $ghgh=gghh$. Multiplying with $ g^{-1}$ on the left and with $h^{-1}$ on the right, that is, plugging this equation into the term $g^{-1} x h^{-1}$ and using rule \eqref{rule5} again, shows $hg=gh$.
\end{ex}

\begin{rem}
Notice that an equational proof does not talk about elements of a set. It is a purely \emph{syntactical} deduction. Even when we write ``For all elements $x,y$, ...'' (for the sake of readability of the proof), these are actually just variables. Alternatively, we may interpret an equational proof as a classical proof where every statement has the form ``For all elements ..., the equation ... holds'', and only a few basic deduction rules are allowed, namely the ones in \cref{eqlogic}.

In the case of the algebraic theory of rings, we are not allowed to use ideals, subrings, quotient rings, direct products of rings, let alone more advanced techniques from commutative or non-commutative algebra. The LEM is not allowed either, so a commutativity proof by contradiction (``If the ring is not commutative, then ...'') is not possible. The absence of LEM also entails that the axiom of choice is not available (Diaconescu's Theorem). But in our equational proofs we allow every concept that can easily be reduced to a ``pure'' equational proof. For rings, this includes the arithmetic of integers, modulo arithmetic, and calculations with polynomials. The statement ``$x$ is central'' abbreviates $xy=yx$ for a (free) variable $y$.

Therefore, equational proofs can be seen as the most direct and elementary proofs possible, even though they can be tough to find.  
\end{rem}

\begin{ex} \label{2case}
Every boolean ring is commutative by an equational proof. We start with
\[-x=(-x)^2=x^2=x.\]
Then we calculate
\[x+y=(x+y)^2 = x^2+xy+yx+y^2=x+xy+yx+y\]
and hence $xy=-yx=yx$.

In contrast, the following proof of the same fact is \emph{not} equational: Every reduced ring is a subdirect product of rings without zero-divisors \cite{Kl80}. The only boolean rings without zero divisors are $0$ and $\IF_2$, which are commutative. So the subdirect product must be commutative as well. (Incidentally, this approach reduces Jacobson's theorem to the case of division rings.)

The second proof shows that every \emph{model} of the theory of boolean rings is commutative. That this implies the existence of an equational proof is the content of Birkhoff's Completeness Theorem \cite{Bi35}. But as mentioned in the introduction, this approach does not allow us to \emph{find} an equational proof. From a constructive point of view, it is pretty much useless.
\end{ex}
 
\begin{defi}
Let $A$ be a rng.
\begin{enumerate}
\item We call $A$ \emph{potent} when for every $x \in A$ there is some natural number $n>1$ with $x^n=x$.
\item Let $n > 1$. We call $A$ an \emph{$n$-rng} (or \emph{$n$-ring} in the unital case) when $x^n=x$ holds for every $x \in A$.
\end{enumerate}
\end{defi}

\begin{rem}
\indent
\begin{enumerate}
\item The class of $n$-rings has been studied for example in \cite{MM37}, \cite[Appendix B]{Sch10} (under the name \emph{$n$-boolean rings}) and \cite[Section I.12]{Pi67}. If $p \in \IP$, some authors require $pA=0$ in the definition of a $p$-ring, but we will not do that.
\item There is an algebraic theory $\T_n$ whose models are $n$-rings: just add the axiom $x^n=x$ to the ring axioms. In particular, \cref{eqlogic} applies and we can speak of equational proofs in the theory of $n$-rings. 
\item The definition of an $n$-ring only uses the multiplicative structure, so that \emph{$n$-monoids} can be defined in the obvious way. But they are not automatically commutative. The smallest counterexample is the unitalization of the semigroup $\{a,b\}$ with $x \cdot y \coloneqq x$.
\item Recently, Oman \cite{Om23} has shown that a rng is potent if and only if every non-zero subrng contains a non-zero idempotent.
\item Jacobson's Theorem states that every potent rng, and hence every $n$-rng, is commutative \cite{Ja45}. We refer to \cite{Wa71,NT74,Do76} for elementary (albeit non-equational) proofs. We will see in a moment that it is actually sufficient to consider the unital case.
\end{enumerate}
\end{rem}
 
\begin{rem}
Let us gather some basic observations on potent rings and $n$-rings.
\begin{enumerate}
\item For every prime power $q$ the finite field $\IF_q$ is a $q$-ring. This follows from Lagrange's Theorem applied to $\smash{\IF_q^{\times}}$.
\item Every potent rng is reduced. This is because $x^2 = 0 \implies x=0$ holds.
\item In particular, there is an equational proof that every $n$-rng is reduced (meaning that for every term $\tau$ an equational proof of $\tau^2=0$ leads to an equational proof of $\tau=0$).
\item The class of $n$-rings is closed under product rings, quotient rings and subrings.
\item The class of potent rings is closed under finite products, but not under infinite products. For instance, $\smash{\prod_{p \in \IP} \IF_p}$ is not a potent ring. One can even show that $(\IF_p)_{p \in \IP}$ has no product in the category of potent rings, which also shows that potent rings are not modelled by an algebraic theory.
\item Every potent ring that is an integral domain is a field. It follows that in a commutative potent ring every prime ideal is maximal, i.e.\ it has Krull dimension $\leq 0$. (Recall that the zero ring has Krull dimension $-\infty$.)
\item When $n$ is even, there is an equational proof that every $n$-rng satisfies $-x=x$.
\end{enumerate}
\end{rem}
 
\begin{lemma} \label{idemcentral}
In a reduced rng, idempotents are central. In an $n$-rng, an equational proof of this claim is available.
\end{lemma}

\noindent Notice that we cannot claim an equational proof in the case of reduced rings, since they do not form an algebraic theory.

\begin{proof}
Let $e^2 = e$ and $x$ be any element. We compute
\[(exe - ex)^2 = exeexe - exeex - exexe + exex =exexe - exex - exexe + exex = 0.\]
Since the rng is reduced, we get $exe = ex$. A similar calculation shows $exe = xe$. Hence, $e$ commutes with $x$. In an $n$-rng, we can make this argument equational: From $(exe-ex)^2=0$ we conclude $exe-ex = (exe-ex)^n = 0$ (since $n > 1$), and proceed exactly as before.
\end{proof}

\begin{lemma} \label{idempower}
Let $A$ be a rng, $n>1$ and $x \in A$ be an element with $x^n=x$. If $k,k' \geq 1$ satisfy $k \equiv k' \bmod n-1$, then $x^k = x^{k'}$ by an equational proof. In particular, $x^{n-1}$ is idempotent by an equational proof.
\end{lemma}
 
\begin{proof}
Based on $x^{n-1} x = x$, a simple induction shows $x^{(n-1)i + j} = x^j$ for all $ i \geq 0$, $j \geq 1$.
\end{proof}

\begin{cor} \label{inclusion}
If $n,m > 1$ are such that $m-1 \mid n-1$, then every $m$-ring is an $n$-ring by an equational proof. \hfill $\square$
\end{cor}

The last three results will be used a lot in the following. Let us give two sample applications:

\begin{ex}\label{3case}
Here is an equational proof that every $3$-rng is commutative, taken from \cite{BM13}. Let $x$ be any element. By \cref{idemcentral,idempower} every square $y^2$ is central. In particular, $(x^2+x)^2$ is central. It simplifies to $x^4+2x^3+x^2=x^2+2x+x^2=2x^2+2x$. Since $x^2$ is central, it follows that $2x$ is central. Applying $(x^2+x)^2=2(x^2+x)$ twice, we get $x^2+x = (x^2+x)^3 = 4(x^2+x)$, hence $3x^2+3x=0$. Since $x^2$ is central, it follows that $3x$ is central. Thus, $x = 3x-2x$ is central as well.
\end{ex}

\begin{ex}\label{4case}
The following equational proof that every $4$-rng is commutative is inspired from \cite{Wa99}. Let $x$ be any element. Then $-x=x$, since $4$ is even. Also, every cube $x^3$ is central by \cref{idemcentral,idempower}. We claim that $f(x) \coloneqq x + x^2$ is idempotent, and hence central. In fact, it squares to $x^2 + 2x^3 + x^4 = x^2 + 0 + x = x + x^2$. For all elements $x,y$ then also $f(x+y)-f(x)-f(y)$ is central, which simplifies to $xy+yx$. In particular, $x(xy+yx)=(xy+yx)x$, which simplifies to $x^2 y = y x^2$. So $y$ commutes with $x^2$. Since $x + x^2$ is central, $y$ also commutes with $x$. (We will generalize this method in \cref{p2case}.)
\end{ex}

\begin{rem} \label{3caselong}
As explained in the introduction, every equational proof of the commutativity of $n$-rings leads to a representation of the commutator $XY-YX$ as a sum of expressions of the form $g \cdot (f^n - f) \cdot h$ in the free ring $\IZ\langle X,Y \rangle$. By unwinding what happens in \cref{3case}, we get the following equation, where $p(f) \coloneqq f^3-f$.
{\small \begin{align*}
& XY-YX = \\
& + \bigl(p(X+X^2) - p(X) (4+4 X+3 X^2+X^3)\bigr) Y  - Y \bigl(p(X+X^2) - p(X) (4+4 X+3 X^2+X^3)\bigr) \\
& + 3 p\bigl(X^2 Y X^2-Y X^2\bigr)  - 3 p\bigl(X^2 Y X^2-X^2 Y\bigr) \\
& + 3 \bigl(X^2 Y X^2-X^2 Y\bigr) X^2 Y X p(X) Y (X^2-1)  - 3 \bigl(X^2 Y X^2-Y X^2\bigr) (X^2-1) Y X p(X) Y X^2 \\
& + p\bigl((X+X^2)^2 Y (X+X^2)^2 - Y (X+X^2)^2 \bigr)  - p\bigl((X+X^2)^2 Y (X+X^2)^2 - (X+X^2)^2 Y \bigr) \\
& + \bigl((X+X^2)^2 Y (X+X^2)^2 - (X+X^2)^2 Y\bigr)   (X+X^2)^2 Y (X+X^2) p(X+X^2) Y \bigl((X+X^2)^2-1\bigr) \\
& - \bigl((X+X^2)^2 Y (X+X^2)^2 - Y (X+X^2)^2\bigr)   \bigl((X+X^2)^2-1\bigr) Y (X+X^2) p(X+X^2) Y (X+X^2)^2 \\
& + 2 p\bigl(X^2 Y X^2-X^2 Y\bigr) - 2 p\bigl(X^2 Y X^2-Y X^2\bigr) \\
& + 2 \bigl(X^2 Y X^2-Y X^2\bigr) (X^2-1) Y X p(X) Y X^2  - 2 \bigl(X^2 Y X^2-X^2 Y\bigr) X^2 Y X p(X) Y (X^2-1) \\
& + (2 + X) p(X) Y  - Y (2 + X) p(X)
\end{align*}}\noindent
It shows that every $3$-ring is commutative. Even though this is a very bad way of writing down the proof, it is interesting that such a single equation always exists. And it also shows for $d \in \IZ$ that a ring satisfying $d x= d x^3$ for all $x$ satisfies $dxy=dyx$ for all $x,y$.
\end{rem}

\begin{prop} \label{unitalred}
Commutativity can be reduced to the unital case. More precisely:
\begin{enumerate}
\item If every potent ring is commutative, then every potent rng is commutative. 
\item If every $n$-ring is commutative by an equational proof, then every $n$-rng is commutative by an equational proof.
\end{enumerate}
\end{prop}

\begin{proof}
1. Let $A$ be a potent rng. Let $x,y \in A$ and choose $n,m>1$ with $x^n=x$ and $y^m=y$. By \cref{idemcentral,idempower} the elements $x^{n-1}$ and $y^{m-1}$ are central. Consider the element
\[e \coloneqq x^{n-1} + y^{m-1} - x^{n-1} y^{m-1} = y^{n-1} + x^{m-1} - y^{n-1} x^{m-1}.\]
We calculate
\[e x = x^{n-1} x + y^{m-1} x - x^{n-1} y^{m-1} x = x^n + y^{m-1} x - x^n y^{m-1} = x + y^{m-1} x - x y^{m-1} = x.\]
A similar calculation shows $e y = y$. By symmetry we also get $x e = x$ and $y e = y$. Hence, if $B$ denotes the subrng of $A$ generated by $x,y$, which consists of $\IZ$-linear combinations of finite non-empty products of $x$ and $y$, then $e$ is multiplicative identity of $B$, and clearly $B$ is potent. By assumption, $B$ is commutative, so that $xy=yx$.

2. In the case of $n$-rngs, we can take $m=n$ and define $e$ as above. There is an equational proof of $xy=yx$ in the theory of $n$-rings. We may assume that $x,y$ are the only variables in the proof. Since $e$ is neutral for $x$ and $y$, every appearance of the multiplicative identity can be replaced by $e$. This yields an equational proof of $xy=yx$ in the theory of $n$-rngs.
\end{proof}

Let us demonstrate this reduction in another special case:

\begin{lemma} \label{binomrel}
In an $n$-rng, there is an equational proof of
\[\sum_{k=1}^{n-1} \binom{n}{k} x^k = 0.\]
\end{lemma}

\begin{proof}
If we had a multiplicative identity, then this follows simply by expanding $1+x=(1+x)^n$ with the binomial theorem and then cancelling $1$ and $x^n=x$. If not, we just do the same with $x^{n-1} + x  = (x^{n-1} + x)^n$ and use that $x^{n-1}$ is neutral for all powers of $x$.
\end{proof}

Because of \cref{unitalred} we will only work with rings (i.e.\ unital rings) from now on.

\begin{rem} \label{unitchar}
The following are equivalent for an element $x$ of an $n$-ring.
\begin{enumerate}
\item $x$ is a unit.
\item $x$ is no zero-divisor.
\item $x^{n-1} = 1$.
\end{enumerate}
\end{rem}

\begin{lemma} \label{unitdecomp}
Let $A$ be an $n$-ring. Then every element $x \in A$ can be written as $x = e \cdot u$, where $e \in A$ is idempotent and $u \in A$ is a unit. Namely, we have $e = x^{n-1}$ and $u = x + (1-x^{n-1})$.
\end{lemma}

\begin{proof}
The element $e \coloneqq x^{n-1}$ is idempotent by Lemma \ref{idempower}. Let $u \coloneqq x + (1-e)$. Then we have $eu = ex =x^n = x$. We need to prove that $u$ is a unit. We use \cref{unitchar}. Notice that the product of $x$ and $1-e$ is zero, and $1-e$ is idempotent. Hence, the binomial theorem shows $u^{n-1} = x^{n-1} + (1-e)^{n-1} = x^{n-1} + (1-e) = 1$. Alternatively, one may verify directly that $x^{n-2} + (1-e)$ is inverse to $u$.
\end{proof}

\begin{cor} \label{unitsuffice}
Let $A$ be an $n$-ring. If every unit of $A$ commutes with every unit of $A$ (by an equational proof), then $A$ is commutative (by an equational proof).
\end{cor}

\begin{proof}
This follows from \cref{unitdecomp,idemcentral}.
\end{proof}

% Structure of n-rings

\section{Structure of \texorpdfstring{$n$}{n}-rings} \label{sec:struct}

In this section we will investigate the well-known structure of (commutative) $n$-rings. We also find an equational proof of the reduction to prime characteristic.

\begin{rem} \label{nfields}
The $n$-fields are those $\IF_q$, where $q$ is a prime power satisfying $q-1 \mid n-1$. This is because $\smash{\IF_q^{\times}}$ is cyclic of order $q-1$. See Table \ref{fieldtable} for a list of examples. It has been generated programmatically with a bit of SageMath\footnote{\href{https://www.sagemath.org/}{https://www.sagemath.org/}} code, see \cref{app:powers}. It follows in particular that every $n$-ring of characteristic $p$ satisfies $p-1 \mid n-1$ (since it contains $\IF_p$ which is then also an $n$-ring).
\end{rem}
\begin{center}
\begin{table}
\begin{tabular}[c]{ll|ll|ll}
$n$ & $n$-fields & $n$ & $n$-fields & $n$ & $n$-fields \\ \hline
$2$ & $\IF_{2}$ &
$18$ & $\IF_{2}$ &
$34$ & $\IF_{2}, \IF_{4}$ \\

$3$ & $\IF_{2}, \IF_{3}$ &
$19$ & $\IF_{2}, \IF_{3}, \IF_{4}, \IF_{7}, \IF_{19}$ &
$35$ & $\IF_{2}, \IF_{3}$ \\

$4$ & $\IF_{2}, \IF_{4}$ &
$20$ & $\IF_{2}$ &
$36$ & $\IF_{2}, \IF_{8}$ \\

$5$ & $\IF_{2}, \IF_{3}, \IF_{5}$ &
$21$ & $\IF_{2}, \IF_{3}, \IF_{5}, \IF_{11}$ &
$37$ & {\small $\IF_{2},\!\IF_{3},\!\IF_{4},\!\IF_{5},\!\IF_{7},\!\IF_{13},\!\IF_{19},\!\IF_{37}$} \\

$6$ & $\IF_{2}$ &
$22$ & $\IF_{2}, \IF_{4}, \IF_{8}$ &
$38$ & $\IF_{2}$ \\

$7$ & $\IF_{2}, \IF_{3}, \IF_{4}, \IF_{7}$ &
$23$ & $\IF_{2}, \IF_{3}, \IF_{23}$ &
$39$ & $\IF_{2}, \IF_{3}$ \\

$8$ & $\IF_{2}, \IF_{8}$ &
$24$ & $\IF_{2}$ &
$40$ & $\IF_{2}, \IF_{4}$ \\

$9$ & $\IF_{2}, \IF_{3}, \IF_{5}, \IF_{9}$ &
$25$ & {\small $\IF_{2},\!\IF_{3},\!\IF_{4},\!\IF_{5},\!\IF_{7},\!\IF_{9},\!\IF_{13},\!\IF_{25}$} &
$41$ & {\small $\IF_2,\!\IF_3,\!\IF_5,\!\IF_9,\!\IF_{11},\!\IF_{41}$} \\

$10$ & $\IF_{2}, \IF_{4}$ &
$26$ & $\IF_{2}$ &
$42$ & $\IF_2$ \\

$11$ & $\IF_{2}, \IF_{3}, \IF_{11}$ &
$27$ & $\IF_{2}, \IF_{3}, \IF_{27}$ &
$43$ & {\small $\IF_2,\!\IF_3,\!\IF_4,\!\IF_7,\!\IF_8,\!\IF_{43}$} \\

$12$ & $\IF_{2}$ &
$28$ & $\IF_{2}, \IF_{4}$ &
$44$ & $\IF_2$ \\

$13$ & {\small $\IF_{2},\!\IF_{3},\!\IF_{4},\!\IF_{5},\!\IF_{7},\!\IF_{13}$} &
$29$ & $\IF_{2}, \IF_{3}, \IF_{5}, \IF_{8}, \IF_{29}$ &
$45$ & $\IF_2,\IF_3,\IF_5,\IF_{23}$ \\

$14$ & $\IF_{2}$ &
$30$ & $\IF_{2}$ &
$46$ & $\IF_2,\IF_4,\IF_{16}$ \\

$15$ & $\IF_{2}, \IF_{3}, \IF_{8}$ &
$31$ & {\small $\IF_{2}, \IF_{3}, \IF_{4}, \IF_{7}, \IF_{11}, \IF_{16}, \IF_{31}$} &
$47$ & $\IF_2,\IF_3,\IF_{47}$ \\

$16$ & $\IF_{2}, \IF_{4}, \IF_{16}$ &
$32$ & $\IF_{2}, \IF_{32}$ &
$48$ & $\IF_2$ \\

$17$ & $\IF_{2}, \IF_{3}, \IF_{5}, \IF_{9}, \IF_{17}$ &
$33$ & $\IF_{2}, \IF_{3}, \IF_{5}, \IF_{9}, \IF_{17}$ &
$49$ & $\IF_2, \IF_3, \IF_4, \IF_5, \IF_7, \IF_9,$\\

&&&&& $\IF_{13}, \IF_{17}, \IF_{25}, \IF_{49}$ 
\end{tabular}
\caption{List of $n$-fields for $2 \leq n \leq 49$} \label{fieldtable}
\end{table}
\end{center}
\vspace{-5mm}

It turns out that $n$-fields generate all commutative $n$-rings in a suitable sense:
 
\begin{prop} \label{nchar}
A commutative ring is an $n$-ring if and only if it is a subring of a direct product of fields of the form $\IF_q$, where $q - 1 \mid n-1$.
\end{prop}

\begin{proof}
The direction $\Leftarrow$ is clear. To prove $\Rightarrow$, let $A$ be a commutative $n$-ring. Then $A$ is reduced, which means that the intersection of all prime ideals is zero. Hence, the canonical map $\smash{A \to \prod_{\mathfrak{p}} Q(A/\mathfrak{p})}$ is injective, where $\mathfrak{p}$ runs through all prime ideals of $A$. Since $ Q(A/\mathfrak{p})$ (which is just $A/\mathfrak{p}$, since $\dim(A) \leq 0$) is clearly an $n$-field, the claim follows.
\end{proof}

Practically, this means that the class of commutative $n$-rings is already determined by the class of $n$-fields. For example, we have the following improvement of \cref{inclusion}:
 
\begin{cor} \label{nequal}
Let $n,n' > 1$. Every commutative $n$-ring is a commutative $n'$-ring if and only if $q-1 \mid n-1 \implies q-1 \mid n'-1$ holds for all prime powers $q$.
\end{cor}
 
\begin{proof}
The direction $\Leftarrow$ follows from \cref{nchar}. The direction $\Rightarrow$ follows by considering the finite field $\IF_q$.
\end{proof}
 
\begin{ex}
Every $10$-ring is a $4$-ring, and every $14$-ring is a $2$-ring; see also Table \ref{fieldtable}.
\end{ex}

\begin{cor}\label{pow}
Every commutative $n$-ring of prime characteristic $p$ is actually a $p^k$-ring for some $k \geq 1 $. Namely, $k$ is the $\lcm$ of all $d \geq 1$ with $p^d-1 \mid n-1$.
\end{cor}
  
\begin{proof}
Let $A$ be a commutative $n$-ring of characteristic $p$. We have seen in \cref{nchar} that $A$ embeds into a product of fields of the form $\IF_{p^d}$ with $p^d - 1 \mid n - 1$. But then $d \mid k$ and therefore $\IF_{p^d} \subseteq \IF_{p^k}$ is also a $p^k$-ring.
\end{proof}

\begin{rem}\label{converse}
If $n,p,k$ are as in \cref{pow}, it is not true that, conversely, every $p^k$-ring of characteristic $p$ is an $n$-ring, since we do not necessarily have $p^k - 1 \mid n-1$. The smallest counterexample is $n=22$ and $p=2$. Then $k=\lcm(1,2,3)=6$. Thus, $\IF_{2^6}$ is an example of a $2^6$-ring that is not a $22$-ring.
\end{rem}

We can turn any ring into an $n$-ring as follows.

\begin{defi} \label{univquotient}
If $R$ is a ring, we define
\[R_{[n]} \coloneqq R/ \langle x^n - x : x \in R \rangle.\]
This is clearly an $n$-ring with a surjective homomorphism $R \to R_{[n]}$. The construction is universal: If $A$ is an $n$-ring and $R \to A$ is a homomorphism of rings, then it lifts uniquely to a homomorphism of $n$-rings $R_{[n]} \to A$. Therefore, we may call $R_{[n]}$ the universal $n$-ring quotient of $R$. 
\end{defi}

\begin{ex} \label{charcompute}
Let us compute the universal $n$-ring quotient of $\IZ$. We have
\[\IZ_{[n]} \coloneqq \IZ/\langle z^n - z : z \in \IZ \rangle = \IZ/{\gcd(z^n - z : z \in \IZ)}\IZ.\]
Let $c$ be the $\gcd$ in question. Of course, $c > 0$. Since $\IZ_{[n]}$ is an $n$-ring and hence reduced, $c$ is square-free. Let $p \mid c$ be a prime factor. This means $z^n \equiv z \bmod p$ for all $z \in \IZ$. Equivalently, $\IF_p$ is an $n$-ring, which means $p-1 \mid n-1$ by \cref{nfields}. Therefore, $c$ is the product of all primes $p$ with $p-1 \mid n-1$. It follows that if $A$ is an $n$-ring, then $c=0$ holds in $A$, since the unique homomorphism $\IZ \to A$ lifts to a homomorphism $\IZ/c\IZ = \IZ_{[n]} \to A$. In other words, the characteristic of $A$ divides $c$.
\end{ex}

\begin{cor} \label{decomp}
Every $n$-ring is a finite product of $n$-rings of prime characteristic.
\end{cor}

\begin{proof}
Let $A$ be an $n$-ring. By \cref{charcompute} its characteristic is $p_1 \cdots p_s$ with distinct prime numbers $p_1,\dotsc,p_s$. Then the Chinese Remainder Theorem (which holds for all rings and sequences of pairwise coprime two-sided ideals) tells us that the canonical homomorphism $A = A/p_1 \cdots p_s A  \to A/p_1 A \times \cdots \times A/p_s A$ is an isomorphism.
\end{proof}

As a consequence, the commutativity problem can be reduced to prime characteristic. We need to convince ourselves, though, that there is also an equational proof of this reduction.
 
\begin{rem} \label{charmeaning}
Even though it is tempting to say that every boolean ring has characteristic $2$, this is not true: The trivial ring is a boolean ring of characteristic $1$. Every non-trivial boolean ring has characteristic $2$. But since the LEM is not available for equational proofs, we cannot (and do not need to!) decide if a given ring is trivial or not. There is no need to exclude the trivial ring. Therefore, in the following, we will not speak of rings of characteristic $p$, and instead only add $p=0$ as an axiom (a shorthand for $p \cdot 1_A = 0$). Equivalently, we work with $\IF_p$-algebras. Moreover, notice that in the proof of \cref{decomp} the rings $A/p_i A$ do not necessary have characteristic $p_i$, since they could be trivial. They just satisfy $p_i=0$. Of course, \cref{decomp} remains true (in classical logic).
\end{rem}

\begin{lemma} \label{chareq}
Let $A$ be an $n$-ring. Let $p_1,\dotsc,p_s$ be the prime numbers $p$ with $p-1 \mid n-1$. Then there is an equational proof of $p_1 \cdots p_s =0$ in $A$.
\end{lemma}

\begin{proof}
We saw in \cref{charcompute} that $p_1 \cdots p_s$ is the $\gcd$ of all $z^n -z$ with $z \in \IZ$. The argument used finite fields, which is not an equational concept, but if $n$ is fixed we can ignore this and just compute the $\gcd$ of the numbers $z^n - z$ for $z=1,2,3,\dotsc$ with the Euclidean algorithm until we reach $p_1 \cdots p_s$. It suffices to consider prime numbers for $z$. Then we use the extended Euclidean algorithm to write $p_1 \cdots p_s$ as a $\IZ$-linear combination of the numbers $z^n - z$. These numbers vanish in $A$, so that $p_1 \cdots p_s$ vanishes as well.
\end{proof}
SageMath can do the computation for us for every value of $n$, see \cref{app:char} for the code.
 
\begin{ex} 
Consider $n=7$. The prime numbers $p$ with $p-1 \mid 7-1$ are $p=2,3,7$. We find $\gcd(2^7-2,3^7-3)=42 = 2 \cdot 3 \cdot 7$ and $42 = -17 (2^7-2) + (3^7-3)$. Hence, $2 \cdot 3 \cdot 7 = 0$ holds in any $7$-ring. For $n=3$, the primes are $2,3$, and we already have $2^3-2 = 6 = 2 \cdot 3$.
\end{ex}

\begin{lemma} \label{modproof}
Let $\T$ be an algebraic theory extending ring theory by some axioms, and let $p \in \IZ$. If $\sigma,\tau$ are terms in $\T$ of the same arity such that $\T,\, p = 0 \vdash \sigma=\tau$, then there is a term $u$ of the same arity with $\T \vdash \sigma = \tau + p \cdot u$.
\end{lemma}

\begin{proof}
This is a simple induction on the structure of an equational proof as in \cref{eqlogic}. We only show the most interesting one, that is rule \eqref{rule5}. Assume that $\T$ and $ p=0$ prove the equation $\rho(\sigma_1,\dotsc,\sigma_n) = \rho(\tau_1,\dotsc,\tau_n)$ by means of $\sigma_i=\tau_i$ for all $i$. By induction hypothesis, there are terms $u_1,\dotsc,u_n$ such that $\T$ proves $\sigma_i = \tau_i + p \cdot u_i$. Then $\T$ proves the equation $\rho(\sigma_1,\dotsc,\sigma_n) = \rho(\tau_1 + p \cdot u_1,\dotsc,\tau_n + p \cdot u_n)$. Expanding the right hand side yields a term $v$ with $\rho(\tau_1 + p \cdot u_1,\dotsc,\tau_n + p \cdot u_n) = \rho(\tau_1 ,\dotsc,\tau_n ) + p \cdot v$, as desired. Formally, this requires an induction on the structure of $\rho$.
\end{proof}

\begin{thm} \label{proofdecomp}
Assume that for every $p \in \IP$ with $p-1 \mid n-1$ there is an equational proof that every $n$-ring with $p=0$ is commutative. Then there is an equational proof that every $n$-ring is commutative.
\end{thm}

Because of this result, we may restrict our attention to rings with $p=0$ from now on.

\begin{proof}
Let $A$ be an $n$-ring. Let $p_1,\dotsc,p_s$ be the set of $p \in \IP$ satisfying $p-1 \mid n-1$. Then $p_1 \cdots p_s = 0$ holds in $A$ by \cref{chareq}. If $x,y \in A$, by assumption $xy-yx=0$ can be derived with the additional axiom $p_i=0$. Hence, \cref{modproof} yields $u_i \in A$ with $xy - yx = p_i \cdot u_i$. In classical logic, this step of the proof would be easier: just use that $[x]$ and $[y]$ commute in the quotient ring $A/p_i A$.

Therefore, it suffices to prove the following in equational logic: If $p_1,\dotsc,p_s$ are pairwise coprime integers and $x \in A$, $u_i \in A$ satisfy $x = p_i \cdot u_i$ for $i=1,\dotsc,s$, then we can construct an element $v$ with $x = p_1 \cdots p_s \cdot v$. The case $s=1$ is easy. If $s=2$, the extended Euclidean algorithm yields $ q_1,q_2 \in \IZ$ with $q_1 p_1 + q_2 p_2 = 1$. Then
\[x = q_1 p_1 x + q_2 p_2 x = q_1 p_1 p_2 u_2 + q_2 p_2 p_1 u_1 = p_1 p_2 (q_1 u_2 + q_2 u_1).\]
For the general case, we use induction and know $x = p_1 \cdots p_{s-1} \cdot v$ for some element $v$. Since $p_1 \cdots p_{s-1}$ and $p_s$ are coprime, the case $s=2$ yields an element $w$ with $x = p_1 \cdots p_{s-1} p_s \cdot w$, and we are done.
\end{proof}

\begin{ex} \label{5decomp}
Consider $n = 5$. The primes $p$ with $p-1 \mid 5-1$ are $2,3,5$, and $2^5 - 2 = 2 \cdot 3 \cdot 5$. Hence, $2 \cdot 3 \cdot 5 = 0$ holds in any $5$-ring. Assume that there is a commutativity proof modulo each of these primes. Then we get elements $u_1,u_2,u_3$ with $xy - yx = 2 u_1$, $xy-yx = 3 u_2$, $xy-yx = 5 u_3$. Hence, $xy-yx = 30 (u_3 - u_1 + u_2) = 0$.
\end{ex}

\begin{rem}\label{geochar}
By Stone duality \cite{St36}, the category of boolean rings is anti-equivalent to the category of totally disconnected compact Hausdorff spaces (Stone spaces). There is a similar classification for $n$-rings, which is thus much more precise than \cref{nchar}. Namely, if $p^{d_1},\dotsc,p^{d_s}$ are the powers of $p$ with $p^{d_i} - 1 \mid n  - 1$ and $k \coloneqq \lcm(d_1,\dotsc,d_s)$, then the category of $n$-rings with $p=0$ is anti-equivalent to the category of Stone spaces $X$ with a continuous $C_k$-action such that $\smash{X = \bigcup_{i=1}^{s} X^{C_{k/d_i}}}$. The proof of this anti-equivalence will appear elsewhere.
\end{rem}

% Reduction to prime powers
 
\section{Reduction to prime powers}

The goal of this section is to find an equational proof of \cref{pow}.

\begin{thm} \label{reduc}
Let $n>1$ and let $p \in \IP$ with $p-1 \mid n-1$. Let $k$ be the $\lcm$ of all $d \geq 1$ such that $p^{d}-1 \mid n-1$. Then there is an equational proof that every $n$-ring with $p=0$ is a $p^k$-ring. Moreover, if $p^k - 1 \mid n-1$, the converse is also true.
\end{thm}

This is actually a meta-theorem. The proof is a method that enables us to write down, for every fixed pair $(n,p)$, an equational proof that every $n$-ring with $p=0$ is a $p^k$-ring. We will demonstrate this with several examples. Recall from \cref{charmeaning} why we write $p=0$ instead of ``characteristic $p$''. The theorem is a generalization and better explanation of Morita's reduction results \cite{Mo78} mentioned in the introduction. It also includes MacHale's observation in \cite{Ma86} that $2^m+1$-rings with $2=0$ are boolean as a special case.
 
\begin{proof}[Proof of \cref{reduc}]
Since we merely want to show $\smash{x^{p^k}=x}$, it is clear that only polynomials in one variable over $\IF_p$ appear in the proofs, so that commutativity is not an issue here. 
 
The first step is to use the Euclidean algorithm to compute the polynomial
\begin{equation}\label{def:g}
g \coloneqq \gcd(f^n - f : f \in \IF_p[T], \, \deg(f) < n) \in \IF_p[T].
\end{equation}
This computation can be quite cumbersome, but see \cref{perf} below for performance issues. Then one uses the extended Euclidean algorithm to write $g$ as a linear combination of polynomials of the form $f^n - f$. We claim that
\begin{equation} \label{div}
g \mid T^{p^k} - T.
\end{equation}
The following proof is abstract and certainly not equational, but in practice, for every fixed pair $(n,p)$, \cref{div} can just be verified via polynomial division. This why an equational proof still does exist.
 
The universal $n$-ring with $p=0$ on one generator is the quotient $\IF_p[T]_{[n]} = \IF_p[T]/I$, where $I \coloneqq \langle f^n - f : f \in \IF_p[T] \rangle$ (see \cref{univquotient}). We claim that $I$ is generated by those $f^n - f$ with $\deg(f) < n$. In fact, if $f \in \IF_p[T]$ is arbitrary, then polynomial division yields a polynomial $f'$ of degree $<n$ with $f \equiv f' \bmod T^n-T$, which also implies $f^n-f \equiv f'^n-f' \bmod T^n-T$. Since $\IF_p[T]/I$ is a $p^k$-ring by \cref{pow}, we have $\smash{T^{p^k}-T \in I = \langle g \rangle}$, which proves $\smash{g \mid T^{p^k}-T}$.

Since we already wrote $g$ as a linear combination of polynomials of the form $f^n-f$, we can now do the same for $\smash{T^{p^k}-T}$. So, we have found an equation of the form
\[T^{p^k}- T = \sum_{i=1}^{s} u_i \cdot (f_i^n - f_i)\]
in $\IF_p[T]$. Now, if $A$ is any ring with $p=0$, then for $x \in A$ we get
\[x^{p^k}- x = \sum_{i=1}^{s} u_i(x) \cdot (f_i(x)^n - f_i(x)).\]
For a pure equational proof that is independent from the proof above, just expand the right hand side of this equation and then simplify. The result will be $\smash{x^{p^k}-x}$. If $A$ is an $n$-ring with $p=0$, the right hand side vanishes, so that $\smash{x^{p^k}=x}$.

Finally, if $p^k - 1 \mid n-1$, then every $p^k$-ring is an $n$-ring by \cref{inclusion}.
\end{proof}

\begin{cor} \label{fullred}
In order to give an equational proof that $n$-rings are commutative, one may restrict to $p^k$-rings with $p=0$.
\end{cor}

\begin{proof}
This follows from \cref{proofdecomp} and \cref{reduc}.
\end{proof}

\begin{rem} \label{perf}
The definition of the $\gcd$ in \cref{def:g} requires all $p^n$ polynomials $f \in \IF_p[T]$ of degree $<n$, which is not practical for computations. But a subset of these polynomials is usually sufficient. Let us explain this in more detail.

First we may assume that $f$ is monic (since $p-1 \mid n-1$, constants $u \in \IF_p$ satisfy $u^n=u$ anyway) and $\deg(f) > 0$. Choose an enumeration $f_1,f_2,f_3,\dotsc$ of all monic polynomials of degree $<n$ with $\deg(f_i) \leq \deg(f_{i+1})$. It is recommended to start with $f_1 = T $ and $f_2 = T + 1$. Then the ``partial $\gcd$''
\[g_i \coloneqq \gcd(f_1^n-f_1,\dotsc,f_i^n-f_i),\]
can be computed recursively via $g_{i+1} = \gcd(g_i,f_{i+1}^n-f_{i+1})$. If some $g_i$ already satisfies $\smash{g_i \mid T^{p^k}-T}$, we are done, even when maybe $g_i \neq g$. If not, compute $g_{i+1}$ and continue. It turns out (see the examples below and \cref{lineargcd}) that in many cases a small $i$ will do the job, and that only a few, mostly linear polynomials are required for the $f_1,\dotsc,f_i$.
  
It would be very much desirable to give a ``uniform'' equational proof of \cref{reduc}, which thus can be written down without any specific choice of $(n,p)$ and also without any tedious calculations. But it is not clear if it exists.

In any case, using a computer algebra system does help to find the equations. For the examples below, we have again used SageMath, the code can be found in \cref{app:reduc}.
\end{rem}
 
\begin{rem} \label{gform}
Because of \cref{converse}, the polynomial $g$ is not always equal to $\smash{T^{p^k}- T}$, even though in the examples below it will often be the case. In general, one can show that $g$ is the product of all monic irreducible polynomials in $\IF_p[T]$ whose degree $d$ satisfies $p^d - 1 \mid n-1$. For $k=1$ this means $g = T^p - T$.
\end{rem}

We start with some examples where linear polynomials are sufficient for computing $g$.
  
\begin{ex} \label{reduc-3-2}
Consider $n=3$ and $p=2$. Then $k=1$. The $\gcd$ of $T^3-T$ and $(T+1)^3 - (T+1)$ in $\IF_2[T]$ is $T^2-T$ with the linear combination
\[T^2-T = (T^3-T) + ((T+1)^3-(T+1)).\]
If $A$ is any ring with $2=0$ and $x \in A$, the equation
\[x^2-x = (x^3-x) + ((x+1)^3-(x+1))\]
can simply be verified by hand. This shows that any $3$-ring with $2=0$ is a $2$-ring. Here is a more compact version: We have $x+1=(x+1)^3=x^3+3x^2+3x+1=x+x^2+x+1$, hence $x^2=x$. Alternatively, \cref{binomrel} gives directly $3x + 3x^2 = 0$, hence $x^2=x$.
\end{ex}
 
\begin{ex} \label{reduc-5-2}
Consider $n=5$ and $p=2$. Then $k=1$ and the extended Euclidean algorithm yields
\[T^2-T = (T+1) \cdot (T^5-T) + T \cdot ((T+1)^5-(T+1)) \in \IF_2[T].\]
This produces the following equational proof that every $5$-ring with $2=0$ is a $2$-ring: We have $x+1=(x+1)^5 = x^5+x^4+x+1$, which gives $x^4=x^5=x$, hence $x=x^5=x^2$. Alternatively, we may argue with \cref{binomrel} again and get $0=5x+10x^2+10x^3+5x^4 = x+x^4$, so that $x^4=x$ and then $x=x^2$.
\end{ex}

\begin{ex} \label{reduc-5-3}
Consider $n = 5$ and $p = 3$. Then $k=1$ and
\[T^3 - T = T \cdot  (T^5-T) - (T+1) \cdot  ((T+1)^5-(T+1)) \in \IF_3[T].\]
This leads to the following equational proof that every $5$-ring with $3=0$ is a $3$-ring: We compute $x+1=(x+1)^5 = x^5 + 2x^4 + x^3 + x^2 + 2x + 1 = 2x^4+x^3+x^2+1$, hence $x^4=x^3+x^2-x$ ($\dagger$). Multiplying with $x$ and using $x^5=x$ gives $x=x^4+x^3-x^2$. Plugging this back into ($\dagger$) gives $x^4=x^2$ and hence $x=x^5=x^3$.
\end{ex}

\begin{ex}
Consider $n=6$, hence $p=2$ and $k=1$. Then
\[T^2-T = T^2 \cdot  (T^6-T) + (T^2+1) \cdot  ((T+1)^6 - (T+1)) \in \IF_2[T]\]
leads to the following equational proof that every $6$-ring is a $2$-ring: Since $6$ is even, we have $-x=x$. Next, we compute $x+1 = (x+1)^6 = x^6 + x^4 + x^2 + 1 = x + x^4 + x^2 + 1$, hence $x^4=x^2$. Thus, $x=x^6 = x^2 x^4 = x^2 x^2 = x^4 = x^2$. Conversely, every $2$-ring is a $6$-ring by \cref{inclusion}.
\end{ex}

\begin{ex}
Consider $n = 9$ and $p = 5$. Then $k = 1$, and
\[T^5-T = 2T \cdot (T^9 - T) - (2T + 2) \cdot ((T+1)^9 - (T+1)) \in \IF_5[T]\]
can be used to give an equational proof that every $9$-ring with $5=0$ is a $5$-ring.
\end{ex}

\begin{ex}
Consider $n=10$, hence $p=2$ and $k=2$. Then
\[T^4-T =  (T^2+1) \cdot (T^{10}-T) + T^2 \cdot ((T+1)^{10}-(T+1)) \in \IF_2[T]\]
shows that every $10$-ring is a $4$-ring. Conversely, every $4$-ring is a $10$-ring by \cref{inclusion}.
\end{ex}

\begin{ex}
Consider $n = 17$ and $p=3$. Then $k=2$ and
\[T^9 - T = T \cdot (T^{17}-T) - (T+1) \cdot ((T+1)^{17}-(T+1)) \in \IF_3[T]\]
shows that every $17$-ring with $3=0$ is a $9$-ring. It produces the following equational proof: \cref{binomrel} together with $\smash{\binom{17}{i} \equiv (-1)^i \bmod 3}$ for $0 \leq i \leq 8$ yields
{\small\begin{equation*}
-x^{16}+x^{15}-x^{14}+x^{13}-x^{12}+x^{11}-x^{10}+x^9+x^8-x^7+x^6-x^5+x^4-x^3+x^2-x=0.
\end{equation*}
}\noindent
Multiplying this with $x$ and using $x^{17}=x$ gives
{\small
\begin{equation*}
+x^{16}-x^{15}+x^{14}-x^{13}+x^{12}-x^{11}+x^{10}+x^9-x^8+x^7-x^6+x^5-x^4+x^3-x^2-x=0.
\end{equation*}
}\noindent
Adding these two equations gives $2x^9-2x=0$, so that $x^9=x$.
\end{ex}

\begin{ex}
Consider $n = 35$. Then $p=3$ and $k=1$. In $\IF_3[T]$ we find that $T^3-T$ is the $\gcd$ of $T^{35}-T$ and $(T+1)^{35}-(T+1)$, so that every $35$-ring is a $3$-ring. The converse also holds by \cref{inclusion}.
\end{ex}

\begin{ex} \label{7char2}
Consider $n=7$ and $p=2$. Then $k=2$. Here, the $\gcd$ with linear polynomials does \emph{not} suffice. But we are lucky with a quadratic one:
\begin{align*}
T^4 - T   & = (T^8+T^6+T^5+T^3+1) \cdot (T^7-T) \\
& \quad + (T^8+T^6+T^5+T^3+1) \cdot  ((T+1)^7-(T+1)) \\
& \quad + ((T^2+T+1)^7-(T^2+T+1))
\end{align*}
There is a simpler expression with a cubic polynomial:
\[T^4 - T = (T^{15}+T^{12} + T^3 + 1) \cdot (T^7-T) + T \cdot ((T^3+1)^7 - (T^3+1))\]
We get the following equational proof that every $7$-ring with $2=0$ is a $4$-ring: We compute
\begin{align*}
x^3+1 & = (x^3+1)^7 \\
& = x^{21} + x^{18} + x^{15} +  x^{12} + x^9 +  x^6 +  x^3 + 1 \\
& = x^3 + x^6 + x^3 +  x^6 + x^3 +  x^6 +  x^3 + 1 \quad \text{(by \eqref{idempower})} \\
& = x^6 + 1,
\end{align*}
so that $x^3=x^6$. Multiplying this with $x$ and using $x^7=x$, we get $x^4=x$.
\end{ex}

\begin{ex}
Consider $n = 13$ and $p = 2$. Then $k=2$. Here, we also need linear and quadratic polynomials: In $\IF_2[T]$, one finds that $T^4-T$ is the $\gcd$ of $T^{13}-T$, $(T+1)^{13} - (T+1)$, $(T^2+T+1)^{13} - (T^2+T+1)$. This shows that every $13$-ring with $2=0$ is a $4$-ring.
\end{ex}

\begin{ex}
Every $31$-ring with $2=0$ is a $16$-ring (and vice versa). In this example even polynomials of degree $4$ are required: $T^{16}-T$ is the $\gcd$ of the six polynomials $f^{31}-f \in \IF_2[T]$, where $f$ belongs to $\{T,\, T+1,\, T^2+T+1,\, T^4+T+1,\, T^4+T^3+1,\, T^4+T^3+T^2+T+1\}$.
\end{ex}

In all the examples so far, the polynomial $g$ was just $\smash{T^{p^k}-T}$. But we already mentioned in \cref{gform} that this does not have to be the case. Here is the smallest counterexample: 
 
\begin{ex}
Let $n = 22$, hence $p=2$ and $k=6$. The $\gcd$ of $T^{22} - T$, $(T^2+T)^{22} - (T^2+T)$ in $\IF_2[T]$ is equal to $g' \coloneqq T^{10}+T^9+T^8+T^3+T^2+T$, and this divides $T^{64}-T$. This yields an equational proof that every $22$-ring with $2=0$ is a $64$-ring. (But this is clear anyway because of $21 \mid 63$, and it is unclear if this helps to prove commutativity of $22$-rings.) Incidentally, we have $g'=g$, because \cref{gform} implies that $g = T (T+1)(T^2+T+1)(T^3+T+1)(T^3+T^2+1)=g'$. The complexity of $22$-rings has been already noted by Morita \cite{Mo78}. There is no $m < 22$ such that every $22$-ring is an $m$-ring, as can be seen by considering $\IF_4 \times \IF_8$.
\end{ex}

\begin{rem} \label{lineargcd}
When $k=1$, i.e.\ $p$ is the only power $p^d$ with $p^d - 1 \mid n-1$, let us call $n$ \emph{simple at} $p$. The previous examples and numerical experiments suggest that then linear polynomials are often enough to compute $g$, i.e.\ that $T^p-T = \gcd((T+u)^n - (T+u) : u \in \IF_p)$. Let us call $n$ \emph{nice} for $p$ in this case, otherwise \emph{unpleasant}. See \cref{app:linear} for the SageMath code to find unpleasant numbers. For $p=2$, there are $5484$ numbers $n \leq 10\,000$ that are simple at $2$, and only $142$ of these (that is $2.6\%$) are unpleasant, the smallest example being $n=74$. For $p=3$, all of the simple numbers $n \leq 100\,000$ are nice, and Peter Müller found that the smallest unpleasant number is $n = 1 + (3^{16}-1)/32 = 1\,345\,211$. For $p=5$, the smallest unpleasant number is $n = 1+(5^8 - 1)/3 = 130\,209$. For other prime numbers, the behavior seems to be similar: most simple numbers are nice.
\end{rem}

% The simple case

\section{The simple case}

In this section, we will show how to produce, for every fixed $n > 1$ with a special property (see below), an equational proof of the commutativity of $n$-rings. As before, we will not be able to write down a single equational proof that works for all $n$. The proofs are not ``uniform'' and have to be written down for each $n$ individually. But the recipe is always the same.

\begin{defi} \label{def:simple}
Let us call $n \in \IN$ \emph{simple} if $n > 1$ and every prime power $q$ with $q-1 \mid n-1$ is actually a prime number. The simple numbers $\leq 100$ are (see \cref{app:simple} for the SageMath code to generate them)
\begin{align*}
&2,3,5,6,11,12,14,18,20,21,23,24,26,30,35,38,39,42,44,45,47,48,51,54,\\
&56,59,60,62,66,68,69,72,74,75,77,80,83,84,86,87,90,93,95,96,98.
\end{align*}
They form the OEIS sequence A366343\footnote{\href{https://oeis.org/A366343}{https://oeis.org/A366343}}. The natural density of the set of simple numbers (if it exists) is approximately
\[\lim_{N \to \infty} \frac{\# \{n \leq N : n \text{ is simple}\}}{N} \approx 0.462118,\]
meaning that almost every second number is simple. By \cref{nfields}, a number $n$ is simple precisely when every $n$-field is a prime field. 
\end{defi}

We start by giving a quick proof of Jacobson's Theorem for simple numbers.
  
\begin{thm} \label{jacsimple}
Let $A$ be a ring such that for every $a \in A$ there is some simple number $n > 1$ such that $a^n=a$. Then every element of $A$ is a $\IZ$-linear combination of idempotent elements. In particular, $A$ is commutative.
\end{thm}
 
\begin{proof}
Let $a \in A$. Choose some simple number $n>1$ with $a^n=a$. Let $c$ be the characteristic of $A$, which satisfies $c>0$ (since $2^m=2$ for some $m>0$). The ring $\IZ[a] = (\IZ/c\IZ)[a]$ is finite, because it is generated by $1,\dotsc,a^{n-1}$ as a $\IZ/c\IZ$-module. By applying the Chinese Remainder Theorem to the finite reduced commutative ring $\IZ[a]$, we can write $\IZ[a]$ as a finite direct product of finite fields. If $\IF_q$ is one of these fields, choose some generator $\zeta$ of its multiplicative group and some simple number $m>1$ such that $\zeta^m=\zeta$. Then $\zeta^{m-1}=1$, so that $q-1 \mid m-1$. Since $m$ is simple, $q$ is a prime number. Therefore, $\IZ[a]$ is actually a finite direct product of prime fields. If $e_1,\dotsc,e_s \in \IZ[a]$ are the corresponding idempotent elements, it follows that $a$ is a $\IZ$-linear combination of $e_1,\dotsc,e_s$. Hence, $a$ is central by \cref{idemcentral}.
\end{proof}

In particular, we have a short proof that for simple $n$ every $n$-ring is commutative. In fact, $n$ is simple if and only if every $n$-ring is additively generated by idempotents. The goal of this section is to make this proof even more constructive. Because of \cref{proofdecomp}, it suffices to look at $n$-rings with $p=0$. We look at an example first.
 
\begin{ex} \label{pow5}
The free $5$-ring with $5=0$ on one generator is $\IF_5[T]/(T^5-T)$. First we will write $[T]$ as a linear combination of idempotent (and hence central) elements. The Chinese Remainder Theorem implies that
\[\IF_5[T]/(T^5-T) \to  (\IF_5)^5, \quad [f] \mapsto \bigl(f(0),f(1),f(2),f(3),f(4)\bigr)\]
is an isomorphism. The image of $[T]$ is $(0,1,2,3,4)$, which is the integer linear combination
$1 \cdot (0,1,0,0,0) + 2 \cdot (0,0,1,0,0) + 3 \cdot (0,0,0,1,0) + 4 \cdot (0,0,0,0,1)$ of the canonical idempotents. Their preimages are the classes of
\[e_i(T) := 1-(T-i)^4 \in \IF_5[T]\]
for $0 \leq i \leq 4$. This is because for $u \in \IF_5$ we have $e_i(u)=1$ for $u=i$ and $e_i(u)=0$ for $u \neq i$. Hence,
\[[T] = 1 \cdot [e_1] + 2 \cdot [e_2] + 3 \cdot [e_3] + 4 \cdot [e_4].\]
 This leads to the following equational proof that $5$-rings with $5=0$ are commutative: Let $x$ be any element. A calculation (using only $5=0$) shows that
\[1 \cdot (1-(x-1)^4) + 2 \cdot (1-(x-2)^4) + 3 \cdot (1-(x-3)^4) + 4\cdot (1-(x-4)^4) = x,\]
or equivalently
\[4 \cdot (x-1)^4 + 3 \cdot (x-2)^4 + 2 \cdot (x-3)^4 + 1 \cdot (x-4)^4 = x.\]
Since each $y^4$ is central by \cref{idemcentral,idempower}, it follows that $x$ is central.
\end{ex}

We will now generalize this example.

\begin{lemma} \label{idemgen}
Let $q$ be a prime power. Let $A$ be an $\IF_q$-algebra and $x \in A$. Then there is an equational proof of
\begin{equation}\label{idemliko}
x = \sum_{u \in \IF_q^\times} u \cdot \bigl(1-(x-u)^{q-1}\bigr).
\end{equation}
\end{lemma}

\noindent The equational proof takes place in the theory of $\IF_q$-algebras. This means that the arithmetic of $\IF_q$ is taken for granted.

\begin{proof}
The abstract proof of \cref{idemliko} works as follows: It is enough to verify the equation
\[X = \sum_{u \in \IF_q^\times} u \cdot \bigl(1-(X-u)^{q-1}\bigr)\]
in the universal example $\IF_q[X]$. Both sides are polynomials of degree $<q$ in $X$, so that it suffices to prove that the corresponding polynomial functions $\IF_q \to \IF_q$ agree. This follows from the observation that for $\alpha \in \IF_q$ we have $(\alpha-u)^{q-1}=\delta_{u,\alpha}$.

The equational proof looks as follows: The case $q=2$ is easy to check. Now let $q>2$. We first prove that for all $j \in \IZ$ we have
\begin{equation}\label{sum}
\sum_{u \in \IF_q^\times}  u^j = \begin{cases} -1, & q-1 \mid j, \\ \phantom{-}0, & \text{else} \end{cases}
\end{equation}
in $\IF_q$. The proof is well-known. Choose a generator $\zeta$ of $\smash{\IF_q^\times}$, and let $s$ be the sum in \cref{sum}. Then $\zeta^j s = s$. If $\zeta^j \neq 1$, we get $s = 0$. Otherwise, we have $q-1 \mid j$ and then each summand of $s$ is $1$, so that $s = q-1 = -1$. This proves \cref{sum}. For $j=1$ we get in particular $\smash{\sum_{u \in \IF_q^\times} u = 0}$ (here, we use $q > 2$).

Now, the binomial theorem and \cref{sum} give us
\begin{align*}
\sum_{u \in \IF_q^\times} u \cdot \bigl(1-(x-u)^{q-1}\bigr) & = \sum_{u \in \IF_q^\times} u  \,-\, \sum_{u \in \IF_q^\times} \, \sum_{i=0}^{q-1} \binom{q-1}{i} \cdot u \cdot (-u)^{q-1-i} \cdot x^i \\
& = - \sum_{i=0}^{q-1} \binom{q-1}{i} (-1)^{q-1-i} \cdot \left(\sum_{u \in \IF_q^\times}  u^{q-i}\right) \cdot x^i \\
& =  -\binom{q-1}{1} (-1)^{q-1-1} (-1) \cdot x = x. \qedhere
\end{align*}
\end{proof}

\begin{thm} \label{pcase}
Let $q$ be a prime power. Let $A$ be an $\IF_q$-algebra and $x \in A$. If $x^q = x$, then $x$ can be explicitly written as an $\IF_q$-linear combination of idempotent elements as in \cref{idemliko}. In particular, if $A$ is a $q$-ring, $A$ is commutative by an equational proof.
\end{thm}

\begin{proof}
We apply \cref{idemgen}. If $u \in \IF_q$, then $(x-u)^q = x^q - u^q = x - u$ since $x,u$ commute, and then $(x-u)^{q-1}$ is idempotent by \cref{idempower}. Then $x$ is central by \cref{idemcentral}.
\end{proof}

\begin{rem}
A different equational proof of the commutativity of $p$-rings with $p=0$ has been found by Forsythe and McCoy \cite[Section 3]{FM46}. It uses the Vandermonde determinant and is similar to our proof for $p^2$-rings below (\cref{p2case}).
\end{rem}
 
\begin{cor}  \label{speceq}
For every simple number $n>1$ there is an equational proof that every $n$-ring is commutative.
\end{cor}

\begin{proof}
By \cref{proofdecomp} it is enough to find an equational proof that $n$-rings with $p=0$ are commutative, where $p \in \IP$ satisfies $p-1 \mid n-1$. By \cref{reduc} and the definition of a simple number every such ring is actually a $p$-ring by an equational proof. Now \cref{pcase} finishes the proof.
\end{proof}

Let us check how the equational proof looks like in some examples.

\begin{ex}
Consider $n=3$. First we will use some non-equational concepts, but we will get rid of them afterwards. Let $A$ be a $3$-ring. Then $0 = 2^3 - 2 = 6 = 2 \cdot 3$ holds in $A$. So,
\[A \cong A/2A \times A/3A\]
by the Chinese Remainder Theorem.
The factor $A/2A$
is a $3$-ring with $2=0$, hence a $2$-ring by \cref{reduc-3-2} and therefore commutative. The factor $A/3A$ is a $3$-ring with $3=0$, 
hence commutative by \cref{pcase}.
So $A$ is a product of two commutative rings, and we are done.

Just by unwinding what happens in the proofs, the equational proof looks as follows: It is sufficient to show $xy \equiv yx \bmod 2$ and $xy \equiv yx \bmod 3$, since then $xy \equiv yx \bmod 6$ by the proof in \cref{proofdecomp} and hence $xy=yx$. In fact, if $xy-yx=2u$ and $xy-yx=3v$, then $xy-yx=6(u-v)=0$. We have
\[x+1 = (x+1)^3 = x^3 + 3x^2 + 3x + 1 = 3x^2 + 4x + 1 \equiv x^2+1 \bmod 2\]
and hence $x^2 \equiv x \bmod 2$ for all $x$. Then (we now repeat \cref{2case})
\[x + y \equiv (x+y)^2 = x^2 + xy + yx + y^2 \equiv x + xy - yx + y \bmod 2\]
and hence $xy \equiv yx \bmod 2$. Now for the modulo $3$ part: We have
\[1 \cdot (1-(x-1)^2) + 2 \cdot (1-(x-2)^2) = -3x^2+10x-6 \equiv x \bmod 3.\]
So it suffices to show that $xy^2 \equiv y^2 x \bmod 3$ for all $x,y$. But $y^2$ is central by \cref{idemcentral} and \cref{idempower} (we don't need to repeat the proofs here).
\end{ex}

Of course, there are simpler proofs of the commutativity of $3$-rings (such as the one in \cref{3case}), but the advantage is that our proof has been produced with a \emph{general method}. It was not necessary to come up with clever ideas, and the same method works for larger numbers as well, as long as they are simple. Let us also mention that the reduction to $p$-rings will often be rather simple because of \cref{lineargcd}.

\begin{ex}
Let us show with our general method that every $5$-ring is commutative. We already saw in \cref{5decomp} that it is enough to prove this modulo the primes $2,3,5$. For the primes $2$ and $3$, \cref{reduc-5-2,reduc-5-3} reduce the problem to $2$- and $3$-rings, which we have covered before. For the prime $5$, we can use \cref{pcase}. Explicitly, the equation
\[1 \cdot (1-(x-1)^4) + 2 \cdot (1-(x-2)^4) + 3 \cdot (1-(x-3)^4) + 4 \cdot (1-(x-4)^4) \equiv x \bmod 5\]
shows that $x$ is central modulo $5$ because each $y^4$ is central by \cref{idemcentral,idempower}.
\end{ex}

\begin{ex}
Let us show with our general method that every $21$-ring is commutative. The primes $p$ with $p-1 \mid 21-1$ are $p=2,3,5,11$, we must have $2 \cdot 3 \cdot 5 \cdot 11 = 0$ in any $21$-ring: This is because $2 \cdot 3 \cdot 5 \cdot 11 = 330 = \gcd(2^{21}-2,3^{21}-3,5^{21}-5)$. Thus, it suffices to show commutativity modulo each of these primes. For $ p \in \{2,3,11\}$ we find
\[\gcd\bigl(T^{21}-T,(T+1)^{21}-(T+1)\bigr) = T^p-T\]
in $\IF_p[T]$, showing that every $21$-ring with $p=0$ is a $p$-ring and hence commutative by \cref{pcase}. In $\IF_5[T]$ we find
\[\gcd\bigl((T+u)^{21}-(T+u) : u \in \IF_4\bigr) = T^5-T,\]
showing that every $21$-ring with $5=0$ is a $5$-ring and hence commutative by \cref{pcase}.
\end{ex}

The reduction method can be applied to non-simple numbers as well. Let us look at some examples. More results of this type are summarized in Table \ref{resulttable}.
 
\begin{ex} \label{7case}
Let us show with our general method that every $7$-ring is commutative, even though $7$ is not simple. The primes $p$ with $p-1 \mid 7-1$ are $p=2,3,7$, and since $2 \cdot 3 \cdot 7 = 42 = \gcd(2^7-2,3^7-3)$, it suffices to prove commutativity modulo each of these primes. We already saw in \cref{7char2} that a $7$-ring with $2=0$ is a $4$-ring, hence commutative by \cref{4case}. In $\IF_3[T]$ we find
\[\gcd\bigl(T^{7}-T,(T+1)^{7}-(T+1)\bigr) = T^3-T,\]
showing that every $7$-ring with $3=0$ is a $3$-ring and hence commutative. Finally, every $7$-ring with $7=0$ is commutative by \cref{pcase}.
\end{ex}

\begin{ex} \label{2023case}
Let us show with our general method that every $2023$-ring is commutative. The primes $p$ with $p-1 \mid 2023-1$ are $p=2,3,7$, and since $2 \cdot 3 \cdot 7 = 42 = \gcd(2^{2023}-2,3^{2023}-3)$, it suffices to prove commutativity modulo each of these primes. In $\IF_2[T]$ we find
\[\gcd\bigl(T^{2023}-T,(T+1)^{2023}-(T+1),(T^2+T+1)^{2023}-(T^2+T+1)\bigr) = T^4 - T,\]
so that every $2023$-ring with $2=0$ is a $4$-ring and hence commutative. In $\IF_3[T]$ we find
\[\gcd\bigl(T^{2023}-T,(T+1)^{2023}-(T+1),(T+2)^{2023}-(T+2)\bigr) = T^3 - T,\]
so that every $2023$-ring with $3=0$ is a $3$-ring and hence commutative. Finally, in $\IF_7[T]$ we find
\[\gcd\bigl(T^{2023}-T,(T+1)^{2023}-(T+1)\bigr) = T^7 - T \in \IF_7[T],\]
so that every $2023$-ring with $7=0$ is a $7$-ring and hence commutative by \cref{pcase}.
\end{ex}

% p^2-rings

\section{Commutativity of \texorpdfstring{$p^2$}{p2}-rings}

Morita \cite{Mo78} has shown by direct and rather long calculations that $p^2$-rings are commutative for $p \in \{2,3,5\}$. He also mentioned that the same method (but with much more effort) works for $p=7$. In this section, we will generalize Morita's approach to show the commutativity of $p^2$-rings with $p=0$. The restriction to $p=0$ (and the fact that $p$ will be arbitrary) allows for a considerable simplification of the calculations.

\begin{lemma} \label{invol}
In a $p^k$-ring with $p=0$ there is an equational proof that for every element $x$ the element $\smash{e \coloneqq x + x^p + \cdots + x^{p^{k-1}}}$ is central.
\end{lemma}

\begin{proof}
Observe that $e^p = e$. Hence, \cref{pcase} tells us that $e$ is an $\IF_p$-linear combination of idempotent and hence central elements.
\end{proof}

\begin{lemma} \label{vandermonde}
Let $V$ be an $\IF_p$-vector space, $U \subseteq V$ be a subspace, and $x_0,\dotsc,x_{p-1} \in V$ such that $\smash{\sum_{i=0}^{p-1} \lambda^i x_i \in U}$ for all $\lambda \in \IF_p$. Then $x_0,\dotsc,x_{p-1} \in U$.
\end{lemma}

\begin{proof}
Consider the Vandermonde matrix $(\lambda^i)_{\lambda \in \IF_p,\, 0 \leq i < p}$ over $\IF_p$. It has an explicit inverse matrix $(u_{i,\lambda})_{0 \leq i < p,\, \lambda \in \IF_p}$, so that
\[\sum_{\lambda \in \IF_p} u_{i,\lambda} \lambda^j = \delta_{i,j}.\]
Hence, for all $0 \leq j < p$ we have
\[U \ni \sum_{\lambda \in \IF_p} u_{j,\lambda} \sum_{i=0}^{p-1} \lambda^i x_i
=  \sum_{i=0}^{p-1} \left(\sum_{\lambda \in \IF_p} u_{j,\lambda}  \lambda^i \right) x_i = \sum_{i=0}^{p-1} \delta_{j,i} x_i = x_j. \qedhere\]
\end{proof}

\begin{defi}
Let $A$ be a ring, $x,y \in A$ and $i,j \in \IN$. Following \cite{FM46}, $[x^i y^j]$ denotes the sum of all monomials in $x,y$ in which $x$ appears $i$ times and $y$ appears $j$ times.
There are $\smash{\binom{i+j}{i}}$ of these monomials. For example, $[x y^2]$ equals $x y^2 + yxy + x y^2$. Notice that
\begin{equation} \label{generalbinomi}
(x+y)^n = \sum_{i=0}^{n} [x^i y^{n-i}]
\end{equation}
for all $n \in \IN$. Of course, if $A$ is commutative, then $\smash{[x^i y^j] = \binom{i+j}{i} x^i y^j}$, and \cref{generalbinomi} is just the binomial theorem.
\end{defi}

\begin{thm} \label{p2case}
There is an equational proof that every $p^2$-ring with $p=0$ is commutative.
\end{thm}

\begin{proof}
In the proof, all variables are quantified over all elements of the ring. By \cref{invol}, the element $f(x) \coloneqq x + x^p$ is central. Hence, $f(x+y)-f(x)-f(y)$ is central as well. By \cref{generalbinomi} this simplifies to $[x y^{p-1}] + \cdots + [x^{p-1} y]$. Substituting $x$ with $\lambda x$ for $\lambda \in \IF_p$ shows that
\[\lambda [x y^{p-1}] + \cdots + \lambda^{p-1} [x^{p-1} y]\]
is central. By applying \cref{vandermonde} to the subspace $Z(A) \subseteq A$, all $[x y^{p-1}],\dotsc,[x^{p-1} y]$ are central. In an equational proof, we cannot write $Z(A) \subseteq A$ or speak of subspaces, but it is clear that the proof of the Lemma still applies: we just have to replace the formula $a \in Z(A)$ by $ab=ba$ for a variable $b$. In particular, $[x y^{p-1}]$ commutes with $y$, meaning
\[y(x y^{p-1} + y x y^{p-2} + \cdots + y^{p-2} x y + y^{p-1} x) = (x y^{p-1} + y x y^{p-2} + \cdots + y^{p-2} x y + y^{p-1} x) y.\]
The $i$th summand on the left side is equal to the $(i+1)$th summand on the right side. Hence,
\[y^p x = x y^p,\]
showing that $y^p$ is central. Since $y + y^p$ is central, $y$ must be central.
\end{proof}

\begin{ex} \label{73case}
There is an equational proof that every $73$-ring is commutative: The set of $p \in \IP$ with $p-1 \mid 73-1$ is $2,3,5,7,13,19,37,73$, so by \cref{proofdecomp} it suffices to prove commutativity modulo these primes. We use \cref{reduc}. The exponent $k$ is $1$ for the primes $7,13,19,37,73$ and $2$ for the primes $2,3,5$. So the claim follows from \cref{pcase} for the first set of primes and from \cref{p2case} for the second set of primes.
\end{ex}

% General case

\section{The general case}

Next, we aim for the general case. Our methods for $p$-rings and $p^2$-rings do not generalize to $p^3$-rings, so that different ideas are required. The (non-equational) proofs of Jacobson's Theorem in \cite{Wa71,NT74} derive it from Wedderburn's Theorem, which states that every finite division ring is commutative.

More precisely, they are structured as follows:

\begin{enumerate}
\item Reduce the general case to the case of finite rings.
\item Reduce the finite case to the case of division rings.
\item For division rings, cite (or prove again) Wedderburn's Theorem.
\end{enumerate}

We will follow a similar approach here, but the proofs need to be made constructive (thus, avoiding LEM) and then transformed into equational logic. The common proofs of the Wedderburn Theorem such as the one by Witt \cite{Wi31} are not constructive, though. They start with a finite non-commutative division ring of minimal cardinality and derive a contradiction. To solve this, we need to formulate the Wedderburn Theorem in a different way.

But first, let us observe the following interesting relationship:

\begin{prop} \label{finite}
In classical logic, the following are equivalent:
\begin{enumerate}
\item Wedderburn's Theorem is true.
\item Finite reduced rings are commutative.
\item For every $n > 1$, finite $n$-rings are commutative.
\end{enumerate}
\end{prop}

\begin{proof}
(1) $\Rightarrow$ (2): It suffices to prove that every finite reduced ring $A$ is a product of division rings. If $A$ is trivial or a division ring, we are done. Otherwise (here we are using LEM) we can choose a non-unit $x \neq 0$. Since $A$ is finite, there are $n,m > 0$ with $x^n = x^{n+m}$. Then $e \coloneqq x^{nm}$ is idempotent $\neq 0,1$ and central by \cref{idemcentral}. Thus, $A \cong eA \times (1-e)A$ and induction on the cardinality of $A$ proves the claim. (2) $\Rightarrow$ (3): This is clear since $n$-rings are reduced. (3) $\Rightarrow$ (1): If $D$ is a finite division ring, say with $n$ elements, then $D$ is an $n$-ring because of Lagrange's Theorem applied to $D^{\times}$.
\end{proof}

The idea is now to formulate a constructive version of the Wedderburn Theorem for very special kinds of finite $n$-rings. We already know that it is enough to consider $p^k$-rings with $p=0$ (\cref{fullred}).

\begin{defi} \label{wedderdef}
Let $p \in \IP$, $k \geq 1$, and $f \in \IF_p[T]$ be a polynomial. The \emph{constructive Wedderburn Theorem} $W_{p,k,f}$ is the following statement: Let $A$ be a $p^k$-ring with $p=0$ and $a,b \in A$ two elements with
\begin{equation} \label{commrel}
ba = f(a) b.
\end{equation}
Then
\[ba = ab.\]
This theorem will be an intermediate step in proving commutativity of $A$. It is more easy than just showing $ba=ab$ directly because of the additional assumption in \cref{commrel}, which roughly says that $a$ and $b$ ``almost'' commute. The statement $W_{p,k,T}$ is trivial.
\end{defi}

\begin{lemma} \label{Wc}
The Wedderburn Theorem implies the constructive Wedderburn Theorem $W_{p,k,f}$ for all $p,k,f$. In particular, $W_{p,k,f}$ holds in classical logic.
\end{lemma}

\begin{proof}
With the notation from \cref{wedderdef}, consider the subalgebra $\IF_p\langle a,b\rangle \subseteq A$ generated by $a,b$. Define $\smash{E \coloneqq \{a^i b^j : 0 \leq i,j < p^k\}}$ and let $M$ be the $\IF_p$-submodule of $\IF_p\langle a,b\rangle$ generated by $E$. Then $\smash{a^{p^k}=a}$, $\smash{b^{p^k} = b}$ and \cref{commrel} imply that $a \cdot E$ and $b \cdot E$ are contained in $M$. This implies $E \cdot E \subseteq M$, so that $M \cdot M \subseteq M$ and therefore $M = \IF_p \langle a,b \rangle$. Hence, $\IF_p\langle a,b \rangle$ is a finite $p^k$-ring. So the claim follows from \cref{finite}.
\end{proof}

Trivially, when all $p^k$-rings with $p=0$ are commutative, $W_{p,k,f}$ holds. We will prove next that the converse also holds. The proof is an adaptation of the proofs of Jacobson's Theorem in \cite{Wa71,NT74}. The difference is that we will not use a proof by contradiction, and indeed the proof is (almost) constructive. We will use it to find an equational proof afterwards.

\begin{thm}\label{generalcase}
Let $p \in \IP$ and $k \geq 1$. Assume that the constructive Wedderburn Theorem $W_{p,k,f}$ holds for all $f \in \IF_p[T]$ of degree $<k$. Then every $p^k$-ring with $p=0$ is commutative.
\end{thm}

\begin{proof}
Let $A$ be a $p^k$-ring with $p=0$, and let $a \in A$. Then $\IF_p[a]$ is a finite reduced commutative ring, hence a finite product of finite fields $\IF_{p^{m_i}}$ by the proof of \cref{finite}.
Hence, there are idempotent elements $e_i \in \IF_p[a]$ such that $\sum_i e_i = 1$, $e_i e_j = \delta_{i,j}$ and such that the $\IF_p$-algebra $e_i \cdot \IF_p[a]$ with multiplicative identity $e_i$ is isomorphic to $\IF_{p^{m_i}}$. They are $ p^k$-rings, so that $m_i \mid k$. Notice that $e_i \cdot \IF_p[a]$ is equal to the $\IF_p$-subalgebra of $e_i \cdot A$ generated by $e_i \cdot a$. Since $a = \sum_i e_i \cdot a$, it suffices to prove that every $e_i \cdot a$ is central in $e_i \cdot A$.

Therefore, we may assume that $\IF_p[a]$ is a finite field, say $\IF_{p^m}$ with $m \mid k$. Consider the $\IF_p$-linear map
\[d : A \to A, \quad x \mapsto ax - xa.\]
We need to show $d=0$, because then $a$ is central. One proves by induction on $n \in \IN$ that
\begin{equation*}
d^n(x) = \sum_{i=0}^{n} (-1)^i \binom{n}{i} \cdot a^{n-i} x a^i.
\end{equation*}
Apply this to $n=p^m$ and use $\binom{p^m}{i} \equiv 0 \bmod p$ for $0 < i < p^m$ and $\smash{a^{p^m}=a}$. Then we get
\begin{equation} \label{derivpower}
d^{p^m} = d. 
\end{equation}
Next, consider for every $u \in \IF_p[a]$ the map
\[\lambda_u : A \to A,\quad x \mapsto ux.\]
It is $\IF_p$-linear and commutes with $d$ (since $u$ commutes with $a$). In $\IF_{p^m}[T]$ we have the well-known equation
\[T^{p^m}-T = \prod_{u \in \IF_{p^m}} (T-u).\]
Hence, the same equation holds in $\IF_p[a][T]$:
\begin{equation} \label{powereq}
T^{p^m}-T = \prod_{u \in \IF_p[a]} (T-u).
\end{equation}
Consider the $\IF_p$-algebra homomorphism
\[\alpha : \IF_p[a][T] \to \End_{\IF_p}(A), \quad u \mapsto \lambda_u, \quad T \mapsto d.\]
Here, $\End_{\IF_p}$ refers to the $\IF_p$-algebra of $\IF_p$-module endomorphisms with $\circ$ as the multiplication. Here, $\alpha$ is well-defined since $d$ commutes with each $\lambda_u$. Applying $\alpha$ to \cref{powereq} yields
\[d^{p^m}-d = \prod_{u \in \IF_p[a]} (d-\lambda_u)\]
in $\End_{\IF_p}(A)$. The left side vanishes because of \cref{derivpower}. Therefore, we arrive at
\begin{equation} \label{diffprod}
0 = \prod_{0 \neq u \in \IF_p[a]} (d-\lambda_u) ~ \circ ~ d.
\end{equation}
We claim that $d - \lambda_u : A \to A$ is injective for all $0 \neq u \in \IF_p[a]$. In fact, if $(d - \lambda_u)(b) = 0$, then $ab-ba = ub$, so that $ba = (a-u)b$. Since $a-u \in \IF_p[a]$, there is a polynomial $f \in \IF_p[T]$ of degree $<m \leq k$ with $a-u = f(a)$. Thus, $ba = f(a) b$. Since the constructive Wedderburn Theorem $W_{p,k,f}$ holds by assumption,
it follows $ab=ba$, which means $ub = 0$. Since $u \neq 0$ and $\IF_p[a]$ is a field, $u$ is a unit, so that $b = 0$, and we are done.

Since all $d - \lambda_u$ ($u \neq 0$) are injective, \cref{diffprod} implies $d = 0$, and $a$ is central.
\end{proof}

The next step is to transform the proof into an equational proof.

\begin{thm} \label{generalcase-eq}
Let $p \in \IP$ and $k \geq 1$. Assume that for all $f \in \IF_p[T]$ of degree $<k$ there is an equational proof of the constructive Wedderburn Theorem $W_{p,k,f}$. Then there is an equational proof that every $p^k$-ring with $p=0$ is commutative.
\end{thm}

\begin{proof}
The first step is to describe a product decomposition $\smash{\IF_p[T]/\langle T^{p^k}-T \rangle  \cong \prod_i \IF_{p^{m_i}}}$ in more explicit terms. Let $S$ be the finite set of all monic irreducible polynomials over $\IF_p$ whose degree divides $k$. This can be computed with standard methods. For example, when $p=2$ and $k=2$ we have $S = \{T,\, T+1,\, T^2+T+1\}$. Let $g \in S$. Then $g$ and $\prod_{g' \in S \setminus \{g\}} g'$ are coprime, so that the extended Euclidean algorithm yields polynomials $u,v$ with
\[\textstyle u \cdot g + v \cdot \prod_{g' \in S \setminus \{g\}} g' = 1.\]
Define $e_g \coloneqq v \cdot \prod_{g' \in S \setminus \{g\}} g'$. Then $e_g \equiv 1 \bmod g$ and $e_g \equiv 0 \bmod g'$ for all $g' \in S \setminus \{g\}$. We have
\begin{equation} \label{longprod}
\textstyle T^{p^k}-T = \prod_{g \in S} g.
\end{equation}
The abstract reason is that both sides equal the product of all $T-u$ with $u \in \IF_{p^k}$, but in practice, for fixed $p$ and $k$ we can just verify this by expanding the product. It follows that the images $\smash{[e_g] \in \IF_p[T]/\langle T^{p^k}-T \rangle}$ are idempotent, pairwise orthogonal\footnote{Two idempotents $e,f$ are called orthogonal if $ef = 0$.}, and their sum is $1$.

Let $A$ be a $p^k$-ring with $p=0$, and let $a \in A$. We want to show that $a$ is central. There is a homomorphism of $\IF_p$-algebras
\[\IF_p[T] / \langle T^{p^k}-T \rangle \to A\]
mapping $[T]$ to $a$, and thus $[e_g]$ to $e_g(a)$. It follows that all the elements $e_g(a) \in A$ are idempotent, pairwise orthogonal, and their sum is $1$, meaning
\begin{equation}\label{idemrels}
\textstyle e_g(a)^2 = e_g(a), \quad e_g(a) \cdot e_{g'}(a) = 0 ~ (g \neq g'), \quad \sum_{g \in S} e_g(a) = 1.
\end{equation}
For an equational proof of these equations in $A$, we cannot use homomorphisms and quotients of polynomial rings, but we can just verify them by hand, only using the relation $\smash{a^{p^k}=a}$. In order to show that $\smash{a = \sum_{g \in S} a \cdot e_g(a)}$ is central, it is then enough to show that each $a \cdot e_g(a)$ commutes with each $b \cdot e_g(a)$, where $b \in A$ and $g \in S$. Both elements stay the same when multiplied with $e_g(a)$, which is central by \cref{idemcentral}. Hence, we may assume that $e_g(a)$ is the multiplicative identity of our ring. In fact, we just have to replace every occurrence of $1 \in A$  below by $e_g(a)$, and the calculations remain valid. Of course, this corresponds to the product decomposition in our proof of \cref{generalcase}, but what we do here is an equational version of this argument.

Strictly speaking, we are not allowed to speak of ``irreducible polynomials'' in an equational proof. The construction above, however, just explains the \emph{method} how to write down an equational proof for a fixed pair $(p,k)$, in which we are just using the polynomial expressions $e_g(a)$ and $g(a)$ in $A$. This will become more clear in the examples below.

So, we have a polynomial $g \in S$ with $e_g(a)=1$, and this implies $e_{g'}(a)=0$ for $g' \neq g$. Next, we claim that
\begin{equation} \label{zero}
g(a) = 0
\end{equation}
holds in $A$. To see this, notice that for all $g' \neq g$ we have $1 = e_g(a) \equiv 0 \bmod g'(a)$, i.e.\ $g'(a)$ is a unit. On the other hand, by \cref{longprod} we have $\smash{0 = a^{p^k}-a = \prod_{g' \in S} g'(a)}$. Since all factors $g'(a)$ with $g' \neq g$ are units, we must have $g(a) = 0$. Conversely, $g(a)=0$ implies $e_g(a)=1$ since $e_g \equiv 1 \bmod g$.

The whole preceding discussion was just about finding some $g \in S$  with $g(a)=0$. This is the only thing we will need from it in the following. In classical logic, this means that $\IF_p[a]$ is either zero or a field.

Let $m \coloneqq \deg(g)$. If $m=1$, then $a = u \cdot 1$ for some $u \in \IF_p$ and we are done. So assume $m > 1$ (and hence $k > 1$). Notice that we are \emph{not} using the LEM here. What we are writing down is one proof for every $ g \in S$. When $S = \{T,\, T+1,\, \dotsc\}$, we start with $g = T$, then $g = T + 1$, etc. But for linear polynomials the proof is easy. So we switch directly to the next polynomials. Incidentally, this means that the proof is easy for $k=1$, namely exactly the one in \cref{pcase}.

Since $g$ is irreducible of degree $ > 1$, the constant term of $g$ is a unit in $\IF_p$. So clearly, $a \in A$ is a unit. Since $g$ is irreducible of degree $m$, we have $\smash{g \mid T^{p^m}-T}$ (this can be verified manually with a polynomial division), so that $\smash{a^{p^m}=a}$ by \cref{zero}.

Let $b \in A$. We want to show that $a$ and $b$ commute. Choose an enumeration $f_0,f_1,\dotsc,f_{p^m-1}$ of the set $\IF_p[T]_{<m}$ of polynomials over $\IF_p$ of degree $<m$, starting with $f_0=0$. We define a sequence of elements $b_0,b_1,b_2,\dotsc,b_{p^m}$ in $A$ recursively by
\begin{equation}\label{bndef}
b_0 \coloneqq b, \quad b_{n+1} \coloneqq a b_n - b_n a - f_n(a) b_n.
\end{equation}
For example, $b_1 = a b - b a$, and the goal is to prove $b_1 = 0$. We claim that
\begin{equation} \label{end}
b_{p^m}=0.
\end{equation}
Unfortunately, our proof is not equational, but it is constructive, and \cref{end} can be verified in every example of $p,k,g$ manually as well, just by using $g(a)=0$ (see the examples below). Define the map $d : A \to A$ by $d(x) \coloneqq ax - xa$. As before, we prove $\smash{d^{p^m}=d}$, only using $\smash{a^{p^m}=a}$. The definition of $b_n$ becomes
\[b_n = \bigl((d - f_{n-1}(a)) \circ \dotsc \circ (d - f_0(a))\bigr)(b).\]
Hence,
\[b_{p^m} = \prod_{f \in \IF_p[T]_{<m}} (d - f(a))(b),\]
so that it suffices to prove the equation
\[T^{p^m} - T = \prod_{f \in \IF_p[X]_{<m}} (T - f(a))\]
in $A[T]$. Since $g(a)=0$, there is a homomorphism of $\IF_p$-algebras $\IF_p[X] / \langle g(X) \rangle \to A$, $[X] \mapsto a$, and it suffices to prove
\[T^{p^m} - T = \prod_{f \in \IF_p[X]_{<m}} (T - [f])\]
in $(\IF_p[X]/\langle g(X) \rangle)[T]$. This follows from the fact that $\IF_p[X]/\langle g(X) \rangle$ is a field with $p^m$ elements, namely those $[f]$ with $f \in \IF_p[X]_{<m}$.
This finishes the proof of \cref{end}.

Next, we claim that for $n \geq 1$ there is an equational proof of
\begin{equation} \label{down}
b_{n+1} = 0 \implies b_n = 0.
\end{equation}
If $b_{n+1} = 0$, this means $b_n a = (a - f_n(a)) b_n = f'(a) b_n$ with $f' \coloneqq T - f_n$. Notice that $f'$ has degree $<k$ since $k > 1$ and $f_n$ has degree $<m \leq k$. By assumption, there is an equational proof of the constructive Wedderburn Theorem $W_{p,k,f'}$, so that $b_n a = a b_n$, meaning $f_n(a) b_n = 0$. We claim that $f_n(a)$ is a unit in $A$, which leads to $b_n=0$. In fact, $f_n$ and $g$ are coprime (since $f_n \neq 0$ has degree $ < m = \deg(g)$ and $g$ is irreducible), so that $u f_n + v g = 1$ for some computable polynomials $u,v$. Then $1 = u(a) f_n(a) + v(a) g(a) = u(a) f_n(a)$, proving our claim. 

Combining \cref{end} and \cref{down} gives $b_1 = 0$, so that $a,b$ commute.
\end{proof}

\begin{rem}
If two polynomials $g,g' \in S$ only differ by a linear substitution, then the proof for $g$ in \cref{generalcase-eq} (that $g(a)=0$ implies that $a$ is central)  leads to a proof for $g'$, so that $g'$ may be omitted. In fact, if $g' = w \, g(uT+v)$ for $v \in \IF_p$ and $u,w \in \IF_p^{\times}$, then the assumption $g'(a)=0$ means $g(ua + v)=0$. If there is a proof for $g$, then $ua + v$ is central and hence $a$ too. Similarly, if $g,g'$ are reciprocal polynomials (up to a unit), then the proof for $g$ leads to a proof for $g'$. In fact, if $g'(a)=0$, then $a$ is a unit with $a^m g(a^{-1})=0$, so that $g(a^{-1})=0$. If there is a proof for $g$, then $a^{-1}$ is central, so that $a$ is central as well. Again, one can write some SageMath code that computes a sufficient set of polynomials for us.
\end{rem}

We will now show how the equational proof provided by \cref{generalcase-eq} looks like in the example $p=2$ and $k=2$. It will be considerably more complicated than ours in \cref{p2case}, but it explains the method (that works for \emph{all} $k \geq 1$) quite well.
 
\begin{ex} \label{4casenew}
Consider $p=2$ and $k=2$. So we will show (again) that every $4$-ring $A$ is commutative. We have $S = \{T,\, T+1,\, T^2+T+1\}$. The extended Euclidean algorithm yields the equations
\begin{align*}
T^2 \cdot (T) + 1 \cdot (T+1)(T^2+T+1) & = 1 \\
(T^2+1) \cdot (T+1) + 1 \cdot T(T^2+T+1) & = 1 \\
1 \cdot (T^2+T+1) + 1 \cdot T(T+1) & = 1 
\end{align*}
Therefore, we define
\begin{align*}
e_T & \coloneqq (T+1)(T^2+T+1) = T^3+1,\\
e_{T+1} & \coloneqq T(T^2+T+1) = T^3+T^2+T,\\
e_{T^2+T+1} & \coloneqq T(T+1) = T^2+T.
\end{align*}
Let $a \in A$. Then the elements $e_T(a)$, $e_{T+1}(a)$, $e_{T^2+T+1}(a)$ are idempotent (hence central), pairwise orthogonal, and their sum is $1$. We can verify this by a direct calculation, just using $a^4=a$ (and of course $2=0$ in $A$). We only show two examples: $a^3+1$ is idempotent since $(a^3+1)^2 = a^6+1 = a^3+1$, and $(a^3+1)(a^2+a) = (a^4+a)(a+1)=0(a+1)=0$. It suffices to prove that every $a \cdot e_g(a)$ commutes with all $b \cdot e_g(a)$, where $g \in S$ and $b \in A$. For $g = T$ we find
\[a \cdot e_T(a) = a \cdot (a^3+1) = a^4+a = 0,\]
so this case is trivial. For $g = T+1$ we find
\[a \cdot e_{T+1}(a) = a^4+a^3+a^2 = e_{T+1}(a),\]
so this case is also trivial. Linear polynomials will always be trivial. The only interesting case is $g = T^2+T+1$. Here,
\[a \cdot e_{g}(a) = a^3+a^2.\]
Now, we do not necessarily have $g(a)=a^2+a+1=0$ as in \cref{zero}, but here is why we may assume it: A direct calculation shows
\[(a^3+a^2)^2 + (a^3+a^2) + (a^2+a) = 0.\]
So $c \coloneqq a \cdot e_g(a)$ satisfies $c^2 + c + 1 = 0$ in the $4$-ring $e_g(a) A$ with multiplicative identity $e_g(a)$, and it suffices to prove that $c$ is central in $e_g(e) A$. Formally, this ring does not even exist in equational logic, but any equational proof in $e_g(e) A$ can be turned into an equational proof in $A$ after we replace the multiplicative identity with $e_g(a)$. Alternatively, we may argue as in the proof of \cref{generalcase} and assume $1 = e_{T^2+T+1}(a)= a^2+a$. So we will assume
\[a^2 + a + 1 = 0\]
from now on. We enumerate the polynomials $f_0,f_1,f_2,f_3$ over $\IF_2$ of degree $<2$ by $0,1,T,T+1$. If $b \in A$, we define a sequence recursively by $b_0 \coloneqq b$ and $b_{n+1} \coloneqq a b_n + b_n a + f_n(a) b_n$. We have $b_0 = b$, $b_1 = ab+ba$ and then
\begin{align*}
b_2  & = a b_1 + b_1 a + b_1 \\
& = a(ab+ba) + (ab+ba)a + (ab+ba) \\
& = a^2 b + b a^2 + ab + ba  \quad (\text{using } a^2=a+1) \\
& = (a+1)b + b (a+1) + ab + ba =  0.
\end{align*}
So we are lucky, the sequence terminates quite early. (We only expected $b_4=0$ from the general proof in \cref{generalcase-eq}.) Assume that we can show
\begin{equation*} \label{w2}
ax + xa = x \implies x = 0.
\end{equation*}
Then we can derive $b_1 = 0$ from $b_2 = 0$ and we are done. This is exactly the constructive Wedderburn Theorem $W_{2,2,T+1}$. Let us prove it here: Assume $ax+xa = x$. We calculate
\[(a+x)^2 = a^2 + ax+xa + x^2 = a^2 + (x + x^2).\]
Since $x + x^2$ is idempotent, it is central. Hence, we can calculate
\[a+x = (a+x)^4 = (a^2 + (x + x^2))^2 = a^4 + (x + x^2)^2 = a + x + x^2.\]
From this we conclude $x^2=0$, i.e.\  $x = 0$, and we are done. 

One could remove all references to the more general method from this proof, and also get rid of all remarks outside of equational logic. We will not do that here, also because it will be less clear what is actually happening (and why).
\end{ex}

\begin{ex} \label{8brief}
For $p = 2$ and $k = 3$, one can proceed in a similar way, and we only briefly indicate how it works (we give a better proof in a moment). It suffices to consider $g = T^3 + T + 1$, since the other irreducible polynomial of degree $3$ is $g(T+1)$. So we may assume $a^3+a+1=0$. We enumerate the polynomials $f_i$ of degree $<3$ by $0,\, 1,\, T,\, T+1,\, T^2,\, T^2+1,\, T^2+T,\, T^2+T+1$. Define the sequence $b_n$ recursively by $b_0 \coloneqq b$ and $b_{n+1} \coloneqq a b_n + b_n a + f_n(a) b_n$. Then a routine calculation gives $b_7 = 0$ (even though we only expected $b_8=0$). Hence, if $W_{2,3,f}$ has been verified for all $f$ of degree $<3$, we get an equational proof that $8$-rings are commutative.
\end{ex}

In many examples, one already has $b_n = 0$ for some $n < p^m$. This means that the commutativity proof does not require all $W_{p,k,f}$. Choosing a different enumeration of the $f_i$ will also help to reduce the number $n$ with $b_n=0$. For example, in \cref{8brief} we will get $b_3 = 0$ with a different enumeration. Actually, certain monomials are always sufficient:

\begin{thm} \label{weddermonomial}
Let $p \in \IP$ and $k \geq 1$. Assume that for all $1 \leq i < k$ there is an equational proof of the constructive Wedderburn Theorem $\smash{W_{p,k,T^{p^i}}}$. Then there is an equational proof that every $p^k$-ring with $p=0$ is commutative.
\end{thm}

\begin{proof}
Let $A$ be a $p^k$-ring with $p=0$, and let $a,b \in A$. We will tweak the proof of \cref{generalcase-eq}. Let us first simplify the definition of $b_n$ a bit. The old definition was $b_{n+1} = (a - f_n(a)) b_n - b_n a$. So we might as well define
\begin{equation} \label{bndefnew}
b_0 \coloneqq b, \quad b_{n+1} \coloneqq b_n a - f_n(a) b_n,
\end{equation}
and start the enumeration with $f_0 \coloneqq T$. Again, $W_{p,k,f_n}$ proves $b_{n+1}=0 \implies b_n=0$ for $n \geq 1$ (using that $f_n(a) - a$ is a unit then), and we have $b_1 = b a - a b$. Hence, our goal is to find an equational proof of $b_n=0$ for some $n \geq 1$. As before, we may assume $g(a)=0$ for some monic irreducible polynomial $g \in \IF_p[T]$ whose degree $m$ divides $k$.

Now, even though the computation of $b_n$ takes place inside a ring that is non-commutative a priori, we shall see now that it can be encoded by a computation in a commutative ring. By induction, it follows from \cref{bndefnew} that we can express $b_n$ as
\[\textstyle b_n = \sum_{i,j} \lambda_{i,j} a^i b a^j\]
with $\lambda_{i,j} \in \IF_p$. We encode this with the commutative polynomial
\[\textstyle B_n \coloneqq \sum_{i,j} \lambda_{i,j} X^i Y^j \in \IF_p[X,Y].\]
More formally, we define $B_n \in \IF_p[X,Y]$ recursively by
\[B_0 \coloneqq 1, \quad B_{n+1} \coloneqq B_n Y - f_n(X) B_n .\]
Then $b_n$ results from $B_n$ by replacing any monomial $\lambda X^i Y^j$ with $\lambda a^i b a^j$. This describes a map $\IF_p[X,Y] \mapsto A$ that is clearly additive (but not multiplicative). In particular, $B_n=0$ will imply $b_n=0$. Since $\IF_p[X,Y]$ is commutative, there is a very simple formula for $B_n$: With $B_0 = 1$ and $\smash{B_{n+1} = (Y - f_n(X)) B_n}$ we find
\begin{equation} \label{Bclosed}
B_n = \prod_{i < n} (Y - f_i(X)).
\end{equation}
Consider the finite field $K \coloneqq \IF_p[X]/\langle g(X) \rangle$. (This concept is not equational, but wait for it.) The irreducible polynomial $g(Y) \in \IF_p[Y] \subseteq K[Y]$ has a root in $K$, namely $[X]$. Basic finite field theory tells us that the roots are formed by the orbits under the Frobenius. Thus, we have $\smash{g(Y) = \prod_{i < m} (Y - [X]^{p^i})}$ in $K[Y]$. This means that in $\IF_p[X,Y]$ there is an equation of the form
\begin{equation} \label{gtrans}
g(Y) = \prod_{i < m} (Y - X^{p^i}) +  g(X) \cdot h
\end{equation}
for some $h \in \IF_p[X,Y]$. For every specific choice of $p,g$, we can simply find this equation by calculating $\smash{\prod_{i < m} (Y - X^{p^i})}$, so that no field theory is required. Define $\smash{f_i \coloneqq T^{p^i}}$ for $i < m$; the rest of the polynomials do not matter anymore. Then \cref{Bclosed,gtrans} imply
\[g(Y) = B_m(Y) + g(X) \cdot h.\]
But since $g(a)=0$, we then get $b_m=0$. For a purely equational proof, when $p,g$ are fixed, one can just compute $b_m=0$ directly, just by using $g(a)=0$ (see the examples below).
\end{proof}

\begin{ex} \label{8case}
From \cref{weddermonomial} we get the following equational proof that $8$-rings are commutative, for which only Morita \cite{Mo78} gave an equational proof so far. We may assume $g(a)=0$ for $g \coloneqq T^3+T+1$. For any element $b$, we define a sequence recursively by $b_0 \coloneqq b$ and $\smash{b_{n+1} \coloneqq b_n a + a^{2^n} b_n}$. We calculate:
\begin{align*}
b_1 & = b a + a b, \\
b_2 & = (ba + ab)a + a^2 (ba+ab) \\
& = ba^2 + aba + a^2 ba + a^3 b \quad (\text{use } a^3=a+1) \\
& = b + ab + ba^2 + aba + a^2 ba,\\
b_3 & = (b + ab + ba^2 + aba + a^2 ba)a + (a+a^2) (b + ab + ba^2 + aba + a^2 ba) \\
& = ba + aba + ba^3 + aba^2 + a^2 ba^2 + ab + a^2 b + aba^2 + a^2 ba + a^3 ba \\
& \quad   + a^2 b + a^3 b + a^2 ba^2 + a^3 ba + a^4 ba \\
& = ba + ab + aba + ba^3 + a^2 ba + a^3 b + a^4 ba \quad (\text{use } a^3=1+a, ~ a^4=a+a^2) \\
& = ba + ab + aba + b + ba + a^2 ba + b + ab + aba + a^2 ba \\
& = 0.
\end{align*}
If $W_{2,3,T^4}$ holds, $b_3 = 0$ implies $a^2 b_2 = 0$, hence $b_2 = 0$. If $W_{2,3,T^2}$ holds, then $b_2=0$ implies $(a^2+a) b_1 = 0$, hence $b_1=0$ and $a,b$ commute. We will prove these constructive Wedderburn Theorems in the next section as part of more general results.

We can also replace the proof of $b_3=0$ with a ``commutative computation'', as explained in the general proof in \cref{weddermonomial}. Define polynomials $B_n \in \IF_2[X,Y]$ recursively by $B_0  \coloneqq 1$ and $\smash{B_{n+1} \coloneqq B_n Y + X^{2^n} B_n}$, so that $b_n$ results from $B_n$ by replacing any monomial $X^i Y^j$ with $a^i b a^j$. We compute
\begin{align*}
B_3 & = (Y + X)(Y + X^2)(Y + X^4) \\
    & = Y^3 + (X + X^2 + X^4)Y^2 + (X^3+X^4 + X^6)Y + X^7 \\
		& = Y^3 + X g(X) Y^2 + X^3 g(X) Y +  (X^4+X^2+X+1) g(X) + 1 \\
		& \equiv Y^3 + Y + 1 \bmod g(X) \\
		& \equiv 0 \bmod g(X),g(Y)
\end{align*}
and hence $b_3 = 0$.
\end{ex}

\begin{ex} \label{27case}
With \cref{generalcase-eq,weddermonomial} we get the following equational proof that every $27$-ring with $3=0$ is commutative. (There seems to be no equational proof in the literature yet.) The list of irreducible polynomials over $\IF_3$ whose degree divides $3$ is
\begin{align*}
&T,\, T+1,\, T+2,\, T^3+2T+1,\, T^3+2T+2,\, T^3+T^2+2,\, T^3+T^2+T+2,\\
&T^3+T^2+2T+1,\, T^3+2T^2+1,\, T^3+2T^2+T+1, \,T^3+2T^2+2T+2.
\end{align*}
We may remove the linear polynomials and those polynomials which are linear substitutions of others. This leaves us with $T^3+2T+1$ and $T^3+T^2+2$. So we need to handle two cases.

Assume first that $a^3+2a+1 = 0$. Let us choose an enumeration starting with $T,T+1,T+2$. We calculate:
\begin{align*}
b_1 & = ba - ab,\\
b_2 & = (ba-ab)a - (a+1)(ba-ab) \\
    & = ba^2 + aba + a^2 b - ba + ab,\\
b_3 & = (ba^2 + aba + a^2 b - ba + ab)a - (a+2)(ba^2 + aba + a^2 b - ba + ab) \\
    & = ba^3 + aba^2 + a^2 ba - ba^2 + aba - aba^2 - a^2 ba - a^3 b + aba - a^2 b \\
		& \quad + ba^2 + aba + a^2 b - ba + ab \\
		& = ba^3 - a^3 b - ba + ab \quad (\text{use } a^3=a-1) \\
		& = b(a-1) - (a-1) b - ba + ab = 0.
\end{align*}
Thus, if $W_{3,3,T+1}$ and $W_{3,3,T+2}$ hold, we get $b_1 = 0$, and we are done.
 
Now assume $a^3+a^2+2=0$. Let us choose an enumeration starting with $T,T^3,T^9$. We calculate (we omit the details):
\begin{align*}
b_1 & = ba - ab,\\
b_2 & = b_2 a - a^3 b_2 = b_2  a - (1-a^2) b_2 \\
    & = -b + ab - ba + a^2 b - aba + b a^2 + a^2 ba, \\
b_3 & = b_3 a - a^9 b_3 = b_3 a - (1-a+a^2)b_3  = 0.
\end{align*}
Thus, if $W_{3,3,T^3}$ and $W_{3,3,T^9}$ hold, we get $b_1 = 0$, and we are done.  We will prove the required Wedderburn Theorems in the next section.
\end{ex}

\begin{ex} \label{16case}
Let us prove that $16$-rings are commutative with \cref{weddermonomial}. (Again, there seems to be no equational proof in the literature yet.) The irreducible polynomials over $\IF_2$ whose degree divides $4$ are
\[T^2+T+1,\, T^4+T+1,\, T^4+T^3+1,\, T^4+T^3+T^2+T+1.\]
But $T^4+T+1$ is reciprocal to $T^4+T^3+1$, and $(T+1)^4+(T+1)^3+1 = T^4+T^3+T^2+T+1$. Hence, it is enough to consider the two polynomials $T^2+T+1$ and $T^4+T+1$.

If $a^2+a+1=0$, we proceed as in \cref{4casenew}. Namely, we find that $x \coloneqq ab+ba$ satisfies $xa = (1+a)x$. Assuming $W_{2,4,T+1}$, we get $x=0$ and we are done. Now let $a^4+a+1=0$. If we define $\smash{b_{n+1} \coloneqq b_n a - a^{2^i} b_n}$ for $i=0,\dotsc,4$, then a routine computation, which we omit, shows $b_4 = 0$. Assuming $W_{2,4,2^i}$, we get $b_1 = 0$ and we are done. We will prove the required constructive Wedderburn Theorems in the next section.
\end{ex}

% Constructive Wedderburn Theorems

\section{Constructive Wedderburn Theorems}

In this section, we will study more closely the constructive Wedderburn Theorems $W_{p,k,f}$ from \cref{wedderdef}. Recall that they state that $ba = f(a) b$ holds in a $p^k$-ring with $p=0$ only when $ba=ab$. Because of \cref{generalcase-eq}, it would be very desirable to find equational proofs in general. Here, we will just cover some partial results and special cases. We also indicate how the general case could theoretically be handled by a computer program. Even though we already saw in \cref{weddermonomial} that it is enough to consider $\smash{f = T^{p^m}}$ for $0 \leq m < k$, we will consider arbitrary polynomials $f \in \IF_p[T]$ for the moment.
  
\begin{lemma} \label{wedderconst}
If $f \in \IF_p[T]$ is constant, then there is an equational proof of $W_{p,k,f}$.
\end{lemma}

\begin{proof}
In a $p^k$-ring $A$ with $p=0$, assume $ba=ub$ holds for some $u \in \IF_p$. Then, since $u$ and $b$ commute, we have $b(a-u)=0$. Since $A$ is reduced, this implies $(a-u)b=0$, i.e.\ $ab=ub$. Thus, $ab=ba$.
\end{proof}

For a polynomial $f$ in one variable, we denote by $f^{\circ j}$ the $j$-fold composition of $f$ with itself.

\begin{lemma} \label{powerrel}
If $f \in \IZ[T]$ is a polynomial and the equation $ba = f(a)b$ holds in a ring, then for all $i,j \in \IN$ we have
\[b^{\, j}   a^i = f^{\circ j}(a)^i \, b^{\, j}.\]
More generally, for every $g \in \IZ[T]$ we have
\[b^{\, j}   g(a) = g(f^{\circ j}(a)) \,  b^{\,j}.\]
\end{lemma}

\begin{proof}
The second equation follows from the first. The proof of the first is routine, so a sketch should be enough. One first proves the case $j=1$ by induction on $i$. Then one proves the case $i=1$ by induction on $j$. Then one proves the general case by induction on $i$.
\end{proof}

\begin{lemma} \label{wedderlin}
For $u \in \IF_p$ there is an equational proof of $W_{p,k,T+u}$.
\end{lemma}

\begin{proof}
Assume that $ba = (a+u)b$ holds in a $p^k$-ring with $p=0$. By \cref{powerrel}, we get
\[ba = b^{p^k} a = (a+up^k) b^{p^k} = ab.\qedhere\]
\end{proof}

\begin{lemma} \label{wedderunit}
To give an equational proof of $W_{p,k,f}$, one may assume that $a$ and $b$ are units.
\end{lemma}

\begin{proof}
Let $ba = f(a)b$. Assume that there is an equational proof of $ba=ab$ in case $b$ is a unit. By \cref{unitdecomp} we can write $b = ue$, where $u$ is a unit and $\smash{e = b^{p^{k}-1}}$ is idempotent, hence central. It suffices to prove that $b$ and $ea$ commute, since $b$ commutes with $(1-e)a$ anyway (the products vanish). Consider the $p^k$-ring $eA$ with multiplicative identity $e$ and the homomorphism $\alpha : A \to eA$, $x \mapsto ex$. Then we get $\alpha(b) \alpha(a) = f(\alpha(a)) \alpha(b)$ in $eA$. If $W_{p,k,f}$ has been proven in $eA$, where $b$ is a unit, then $\alpha(b)=b$ commutes with $\alpha(a)=ea$ and we are done. We need to make this argument equational, though, and this works as in the proof of \cref{generalcase-eq}. Observe that $e$ is neutral for both $ea$ and $b$. Any proof of $W_{p,k,f}$ can be applied to the elements $ea$ and $b$ and the multiplicative ``pseudo-identity'' $e$. Thus, $b$ and $ea$ commute, and we are done. The reduction to the case that $a$ is a unit works exactly the same.
\end{proof}

When $b$ is a unit, then $W_{p,k,f}$ states that $ba = f(a)b$ implies $f(a) = a$. Thus, we need to show that $a$ is a fixed point of $f$. We will see in a moment that at least $a$ is a fixed point of a certain power $f^{\circ d}$.

\begin{defi}\label{monoiddef}
Let $f,g \in \IF_p[T]$ and $k \geq 1$. We write $f \equiv_k g$ or just $f \equiv g$ (when $k$ is clear from the context) when there is an equational proof of $f(a)=g(a)$ for all elements $a$ of a $p^k$-ring with $p=0$. The set of equivalence classes becomes a monoid with respect to the composition of polynomials:
\[(f,g) \mapsto f \circ g.\]
Notice that it is well-defined, and the multiplicative identity is $T$. The underlying set can be identified with the one of $\smash{\IF_p[T] / \langle T^{p^k} - T \rangle}$,
since $ f \equiv g$ holds if and only if $\smash{T^{p^k}-T \mid f-g}$. In particular, this monoid is finite, it has exactly $\smash{p^{p^k}}$ elements.
\end{defi}

\begin{rem}
If $f$ is an element of a finite monoid, whose multiplication we write as $\circ$, then there must be a pair of numbers $i \geq 0$ and $\pi > 0$ such that $f^{\circ i} = f^{\circ (i+\pi)}$. Then the powers of $f$ are eventually periodic:
\[1,f,\dotsc,f^{\circ i},f^{\circ (i+1)},\dotsc,f^{\circ (i+\pi-1)},f^{\circ i},f^{\circ (i+1)},\dotsc.\]
We say that $f$ is \emph{of period} $\pi$ and \emph{index} $i$. If $\pi$ is minimal, we call $\pi$ \emph{the} period, similarly for the index.

In particular, every $f \in \IF_p[T]$ has some period and index with respect to the finite monoid $\smash{(\IF_p[T]/\langle T^{p^k}-T \rangle,\circ,T)}$ from \cref{monoiddef}. Both can be efficiently computed for every specific choice of $p,k,f$ because of \cref{idempower}. See also the following example.
\end{rem}

\begin{ex}
\noindent
\begin{enumerate}
\item Constant polynomials are of period $1$ and index $1$.
\item The polynomial $-T$ is of period $2$ and index $0$ (also of period $1$ when $p=2$).
\item For $p=2$, $k=3$, $f = T^2+T$, we compute $f^{\circ 2} = T^4+T$, $f^{\circ 3} = T^4+T^2$ and $f^{\circ 4} = T^8 + T^2 \equiv T^2+T$, so that $f$ is of period $3$ and index $1$.
\item For $p=2$, $k=4$, $f = T^2+1$, we compute $f^{\circ 2} = (T^2+1)^2+1 = T^4$, $f^{\circ 3} = T^8+1$ and $f^{\circ 4}=T^{16} \equiv T$, so that $f$ is of period $4$ and index $0$.
\end{enumerate}
\end{ex}

\begin{prop} \label{periodlemma}
Let $f \in \IF_p[T]$ be a polynomial of period $\pi > 0$. Assume that $a,b$ are two elements of a $p^k$-ring with $p=0$ with $ba = f(a)b$, and $b$ is a unit. Then we have $f^{\circ d}(a)=a$, where $d \coloneqq \gcd(\pi,p^k - 1)$. In particular, if $d=1$, there is an equational proof of $W_{p,k,f}$.
\end{prop}

Notice that this result generalizes both \cref{wedderconst} and \cref{wedderlin}, and that the assumption that $b$ is a unit is harmless because of \cref{wedderunit}. Without this assumption, we can only conclude that $b$ commutes with $a^d$.

\begin{proof}
Choose $i \geq 0$ with $f^{\circ i} \equiv f^{\circ (i+\pi)}$. \cref{powerrel} allows us to compute
\[b^i b^\pi a = b^{i+\pi} a = f^{\circ (i+\pi)}(a) b^{i+\pi} = f^{\circ i}(a) b^i b^\pi = b^i a b^\pi.\]
Since $b^i$ is a unit, this shows that $a$ and $b^\pi$ commute. But then
\[f^{\circ \pi}(a) b^\pi = b^\pi a = a b^\pi.\]
Since $b^\pi$ is a unit, we get
\begin{equation}\label{potj}
f^{\circ \pi}(a) = a.
\end{equation}
Next, we compute
\[f(a) b = ba = b^{p^k} a = f^{\circ p^k}(a) b^{p^k} = f^{\circ p^k}(a) b.\]
Since $b$ is a unit, we get $\smash{f(a) = f^{\circ p^k}(a)}$. Apply $f^{\circ (\pi-1)}$ to both sides and use $a = f^{\circ \pi}(a)$. Then we get
\begin{equation}\label{potp}
f^{\circ (p^k - 1)}(a) = a.
\end{equation}
If $d \coloneqq \gcd(\pi,p^k-1)$, then \eqref{potj}, \eqref{potp} and the B\'{e}zout identity imply $f^{\circ d}(a) = a$. 
\end{proof}

\begin{lemma} \label{monomialperiod}
Let $m \in \IN$ and let $\pi$ be the order of $[m]$ in the additive group $\IZ/k\IZ$. Then the monomial $\smash{T^{p^m}}$ is of period $\pi$ and index $0$.
\end{lemma}

\begin{proof}
Let $\smash{f \coloneqq T^{p^m}}$. Write $\pi m = s k$ with $s \in \IN$. Then $\smash{f^{\circ \pi} = T^{p^{m \pi}} = T^{p^{s k}} \equiv T}$.
\end{proof}

\begin{thm} \label{gcdmain}
Let $p \in \IP$ and $k \geq 1$. If $\gcd(k,p^k-1)=1$, then there is an equational proof that every $p^k$-ring with $p=0$ is commutative.
\end{thm}

\begin{proof}
This follows from \cref{weddermonomial}, \cref{periodlemma} and \cref{monomialperiod}.
\end{proof}

\begin{rem}
For the lack of a better name, let us call $k \geq 1$ \emph{good} for $p \in \IP$ if we have
\[\gcd(k,p^k-1)=1,\]
otherwise \emph{bad} for $p$. With \cref{gcdmain} we have solved the equational commutativity problem for good exponents, so let us analyze how many numbers are good. For $p=2$, the $22$ bad numbers $k \leq 100$ are
\[6, 12, 18, 20, 21, 24, 30, 36, 40, 42, 48, 54, 60, 63, 66, 72, 78, 80, 84, 90, 96, 100.\]
For the other $78$ numbers $k \leq 100$, we thus have an equational proof that $2^k$-rings are commutative. As before, this then leads to many more cases: By \cref{reduc}, every $94$-ring is a $2^{\lcm(1,2,5)}=2^{10}$-ring, hence commutative. For $p > 2$, every good number must be odd. For $p = 3$ and $p = 5$, most odd numbers are also good. For example, the only odd numbers $k \leq 100$ that are bad for $3$ are $39,55$, and the only odd numbers $k \leq 100$ that are bad for $5$ are $55, 93$. But for $p > 5$ there are much more bad odd numbers. If $k \mid k'$ and $k$ is bad, then $k'$ is bad as well, since $\gcd(k,p^k-1) \mid \gcd(k',p^{k'}-1)$.
\end{rem}

The following characterization of bad numbers has been found by Max Alekseyev.

\begin{lemma}
Let $p \in \IP$. The set of bad numbers $k$ for $p$ is equal to the union
\[\bigcup_{q \in \IP,\, q \neq p} q \cdot \ord_{q}(p) \cdot \IN^+.\]
\end{lemma}

\noindent Here, $\ord_q$ refers to the order in the group $(\IZ/q\IZ)^{\times}$.

\begin{proof}
Let $k \geq 1$. If $k$ belongs to the union, there is a prime $q \neq p$ with $q \cdot \ord_{q}(p) \mid k$. Then $q \mid k$ and $q \mid p^k - 1$, so that $\gcd(k,p^k-1) > 1$. Conversely, if $\gcd(k,p^k-1) > 1$, there is a prime $q$ with $q \mid k$ and $q \mid p^k - 1$. The second condition implies $q \neq p$ and then becomes equivalent to $\ord_{q}(p) \mid k$. Since the order divides $q-1$, it is coprime to $q$. Hence, $q \cdot \ord_{q}(p) \mid k$.
\end{proof}

With \cref{gcdmain}, $2^6$ and $3^4$ are the only prime powers $\leq 100$ for which we don't have an equational commutativity proof yet, see also Table \ref{resulttable}. For the moment, the general case seems to be out of reach. But it turns out that $W_{p,k,f}$ can be verified by a computer program, at least theoretically, since it is a \emph{finite} problem:

\begin{prop}
There is a finite $\IF_p$-algebra $U_{p,k,f}$ with a concrete finite presentation such that $U_{p,k,f}$ is commutative if and only if $W_{p,k,f}$ holds.
\end{prop}

\begin{proof}
Because of \cref{wedderunit}, it suffices to work with units. (This is not necessary for the proof, but it will decrease the dimension of the algebra below.) We define the $\IF_p$-algebra by the finite presentation $\smash{V_{p,k,f} \coloneqq  \IF_p\langle X,Y : X^{p^k-1} = 1,~ Y^{p^k-1} = 1,~ YX = f(X) Y \rangle}$. We saw in the proof of \cref{Wc} that $V_{p,k,f}$ is finite, namely $\{X^i Y^j : 0 \leq i,j < p^k - 1\}$ generates $V_{p,k,f}$ as an $\IF_p$-module. Then its universal $p^k$-ring quotient $\smash{U_{p,k,f} \coloneqq (V_{p,k,f})_{[p^k]}}$ from \cref{univquotient} does the job: If two units $a,b$ of a $p^k$-ring with $p=0$ satisfy $ba = f(a)b$, then there is a unique homomorphism of rings $U_{p,k,f} \to A$ that maps $X \mapsto a$ and $Y \mapsto b$. Thus, if $U_{p,k,f}$ commutative, we have $ab=ba$.
\end{proof}

\begin{rem}
With the notation above, $\{X^i Y^j : 0 \leq i,j < p^k-1\}$ is not necessarily an $\IF_p$-basis of $V_{p,k,f}$. For example, $X = 1$ holds in $V_{2,3,T^2}$. More generally, when $p,k$ are arbitrary and $f = T^m$ is any monomial (which is sufficient by \cref{weddermonomial}), then we have $X^d = 1$ for $\smash{d \coloneqq \gcd(m^{p^k-1}-1,p^k-1)}$, and we see that $V_{p,k,f} = \IF_p[G]$ is a group algebra for
\[G \coloneqq \langle X,Y : X^d = 1, \, Y^{p^k-1}=1,\, YXY^{-1} = X^m \rangle \cong C_d \rtimes_m C_{p^k-1}.\]
This description should help to calculate $\smash{U_{p,k,T^m} = \IF_p[G]_{[p^k]}}$, but in practice a brute force attack is still too inefficient, even for small numbers.
\end{rem}

% Appendix with SageMath code

\bigskip

\appendix

\section{SageMath code}

\medskip
\subsection{Computation of \texorpdfstring{$n$}{n}-primes and \texorpdfstring{$n$}{n}-powers} \label{app:powers} 
\phantom{x}\medskip

\noindent The code below has been used to generate the list of $n$-fields in \cref{nfields}.

\medskip
 
\begin{lstlisting}
def n_primes(n):
    """Returns the list of prime numbers p with p-1 | n-1"""
    return [x+1 for x in divisors(n-1) if (x+1).is_prime()]
		
def n_powers(n):
    """Returns the list of prime powers q with q-1 | n-1"""
    return [x+1 for x in divisors(n-1) if (x+1).is_prime_power()]

def all_powers(limit):
    """Returns the lists of n-powers for 2 <= n <= limit"""
    return [[n,n_powers(n)] for n in range(2,limit+1)]
\end{lstlisting}
\medskip

\bigskip
\subsection{Reduction to prime characteristic} \label{app:char}
\phantom{x}\medskip

\noindent The code below implements the algorithm in \cref{chareq}.

\medskip

\begin{lstlisting}
def get_numbers_for_characteristic(n):
    """Returns a list of numbers z such that p_1 * ... * p_s is the gcd of
    the z^n - z, where p_i are the prime numbers with p_i - 1 | n - 1"""
    c = prod(n_primes(n))
    z, g, numbers = 0, 0, []
    while (g != c):
        z = z.next_prime()
        numbers.append(z)
        g = gcd(g, z^n - z)
    return numbers

def characteristic_proof(n):
    """Returns an equational proof that p_1 * ... * p_s = 0 holds in
    any n-ring, where p_i are the prime numbers with p_i - 1 | n - 1"""
    numbers = get_numbers_for_characteristic(n)
    myprimes = n_primes(n)
    relations_string = ", ".join([f"{z}^{n}-{z}" for z in numbers])
    factor_string = " * ".join(map(str,myprimes))
    return f"gcd({relations_string}) = {str(prod(myprimes))} = {factor_string}"
\end{lstlisting}

\medskip

\noindent To get the linear combination, one can use the \begin{tt}extended\_gcd\end{tt} function below.

\bigskip
\subsection{Reduction to prime powers} \label{app:reduc}
\phantom{x}\medskip

\noindent The code below implements the reduction algorithm from \cref{perf}. For example, the command \begin{tt}reduction\_proof(7,2)\end{tt} returns a proof that every $7$-ring with $2=0$ is a $4$-ring.

\medskip

\begin{lstlisting}[language=Python]
def extended_gcd(args):
    """Extension of Sage's xgcd function from two to several arguments.
    Given args = [a_1,....a_n], returns [gcd(a_1,....,a_n),[u_1,....,u_n]]
    where u_1 * a_1 + ... + u_n * a_n = gcd(a_1,...,a_n)."""  
    if (len(args) == 0):
        return [0,[]]
    if (len(args) == 1):
        return [args[0],[1]]
    if (len(args) == 2):
        g, u, v = xgcd(*args)
        return [g,[u,v]]
    first, *rest = args
    rest_gcd, v = extended_gcd(rest)
    full_gcd, u = extended_gcd([first, rest_gcd])
    w = [u[0]] + [u[1] * x for x in v]
    return [full_gcd, w]


def get_exponent(p,n):
    """Returns the least common multiple of all d > 0 with p^d-1 | n-1"""
    limit = floor(log(n,p))
    powers = [d for d in range(1,limit+1) if (p^d-1).divides(n-1)]
    return lcm(powers)


def reduction_proof(n,p,only_polynomials=False):
    """Returns a linear combination of the form g = sum_i u_i * (f_i^n - f_i)
    in GF(p)[T], where g divides T^(p^k) - T, thus proving the reduction theorem.
    When the flag 'only_polynomials' is set, returns only the list of the f_i."""
    if not (p.is_prime() and (p-1).divides(n-1)):
        raise Exception("Invalid arguments")
				
    k = get_exponent(p,n)
    R.<T> = PolynomialRing(GF(p))
    polynomials = []
    g = 0
    target = T^(p^k) - T
    found = False

    # compute gcd of polynomials of the form f^n - f until it divides T^(p^k) - T
    for f in R.polynomials(max_degree = n-1):
        if f.degree() < 1 or (g != 0 and g.divides(f^n-f)):
            continue
        polynomials.append(f)
        g = gcd(g, f^n-f)
        if (g.divides(target)):
            found = True
            break
    if not found:
        raise Exception("No linear combination found")

    if (only_polynomials):
        return polynomials
    
    # print the linear combination
    proof = str(g)
    coeffs = extended_gcd([f^n-f for f in polynomials])[1]
    for i in range(len(polynomials)):
        prefix = " +" if i > 0 else " ="
        coeff = coeffs[i]
        poly = polynomials[i]
        proof += f"{prefix} ({coeff}) * (({poly})^{n} - ({poly}))"
    return proof
\end{lstlisting}

\bigskip
\subsection{Finding unpleasant numbers}  \label{app:linear}
\phantom{x}\medskip

\noindent The code below has been used in \cref{lineargcd} to demonstrate how many simple numbers are nice. We use the \begin{tt}get\_exponent\end{tt} function from above.

\medskip

\begin{lstlisting}[language=Python]
def is_simple_at(n, p):
    """Checks if n is simple at p"""
    return n>1 and (p-1).divides(n-1) and get_exponent(p,n) == 1

def unpleasant_numbers(p,start,end):
    """Returns the number of p-simple numbers within a range and the list of
    unpleasant numbers among them"""
    R.<T> = PolynomialRing(GF(p))
    unpleasants = []    
    simple_count = 0
    for n in range(start,end+1):
        if not (is_simple_at(n,p)):
            continue
        simple_count += 1
        g = gcd([(T+u)^n - (T+u) for u in GF(p)])
        if (g != T^p - T):
            unpleasants.append(n)
    return [simple_count, unpleasants]
\end{lstlisting}

\bigskip
\subsection{Simple numbers} \label{app:simple}
\phantom{x}\medskip

\noindent The code below generates the first simple numbers as in \cref{def:simple}.

\medskip

\begin{lstlisting}
def is_simple(n):
    """Checks if n is a simple number"""
    return n>1 and all(q.is_prime() for q in n_powers(n))

def first_simple_numbers(limit):
   """Returns the list of simple numbers below a given limit"""
   return [n for n in range(2,limit+1) if is_simple(n)]
\end{lstlisting}

\medskip

% References

\end{document}